\documentclass[leqno]{article}
\usepackage{latexsym,amssymb,amsmath,mathrsfs,a4wide}
\usepackage{hyperref}
\newcommand{\p}{\partial}
\newcommand{\f}{\frac}
\newtheorem{theorem}{Theorem}[section]
\newtheorem{lemma}[theorem]{Lemma}
\newtheorem{proposition}[theorem]{Proposition}

\newtheorem{assumption}{Assumption}[section]
\newtheorem{formulation}{Formulation}[section]
\newtheorem{remark}{Remark}[section]

\title{Non-stationary Navier-Stokes Equations with Mixed Boundary Conditions}
\author{Tujin Kim \\
\thanks{Partially supported by AMSS in Chinese Academy of Sciences.}\\
{\small Institute of Mathematics, }\\
{\small Academy of Sciences,}\\
{\small Pyongyang, DPR Korea } \\
\vspace{.5cm}{\small e-mail: math.tujin@star-co.net.kp}\\
Daomin Cao \\
\thanks{Partially supported by NSFC grants (No. 11271354 and No. 11331010) and Beijing
Center for Mathematics and Information Interdisciplinary Sciences.}\\
{\small Institute of Applied Mathematics, }\\
{\small AMSS, Chinese Academy of Sciences,}\\
{\small Beijing 100190, P. R. China} \\
{\small e-mail: dmcao@amss.ac.cn}
 }

\begin{document}

\maketitle

\begin{abstract} In this paper we are concerned with the initial boundary value
problem of the 2, 3-D Navier-Stokes equations with mixed boundary
conditions including conditions for velocity, static pressure,
stress, rotation and Navier slip condition together. Under a
compatibility condition at the initial instance it is proved that
for the small data there exists a unique solution on the given
interval of time. Also, it is proved that if a solution is given,
then there exists a unique solution for small perturbed data
satisfying the compatibility condition.
Our smoothness condition for initial functions in the compatibility condition is weaker than one in such a previous result.
\end{abstract}



\section{Introduction} \setcounter{equation}{0}

For the Navier-Stokes equations
$$
 -\nu\Delta v+(v\cdot\nabla)v+\nabla p=f, \,\,\,
 \nabla \cdot v=0\quad \mbox{in}\,\, \Omega\subset R^l, l=2.3;
$$
and
$$ \f{\p v}{\p t} -\nu\Delta v+(v\cdot\nabla)v+\nabla p=f, \,\,\,
 \nabla \cdot v=0\quad \mbox{in}\,\, \Omega
$$
different natural and artificial boundary conditions are considered.
For example on solid walls, homogeneous Dirichlet condition $v=0$ is
often used. On a free surface a Neumann condition $2\nu
\varepsilon(v)n-pn=0$ may be useful. Here and in what follows
$\varepsilon(v)$ denotes the so-called strain tensor with the
components $\varepsilon_{ij}(v)=\f{1}{2}(\p_{x_i} v_j+\p_{x_j} v_i)$
and $n$ is the outward normal unit vector.  For simulations of flows
in the presence of rough boundaries, the Navier slip-with-friction
boundary conditions $v\cdot
n=0,\,\,(\nu\varepsilon_{n\tau}(v)+\alpha v_\tau)|_{\Gamma_5}=0 $ is
also used, where $\varepsilon_{n\tau}(v)$ and $v_\tau$ are,
respectively, the tangent components of $\varepsilon(v)n$ and $v$.
 Combination of the condition $v_n=0$ and the
tangential component of the friction (slip condition for uncovered
fluid surfaces) or the condition $v_\tau=0$ and the normal component
of the friction (condition for in/out-stream surfaces) are
frequently used. At the outlet of a channel ``do nothing" condition
$\nu \f{\p v}{\p n}-pn=0$, i.e. the outlet boundary condition, is
also used. Rotation boundary condition has been fairly extensively
studied over the past several years.
 Also, for inlet or outlet of flow one deals with the
static pressure $p$ or total pressure (Bernoulli's pressure)
$\f{1}{2} |v|^2+p$. For papers dealing with the problems mentioned
above one can refer to Introduction of \cite{kc}.

In practice we deal with mixture of some kind of boundary
conditions. For a channel flow  a mixture of Dirichlet condition
$v=0$ on the wall and   ``do nothing" condition  on the outlet is
used. But for  a channel flow with a rough boundary surface a
mixture of Dirichlet condition, the Navier slip-with-friction
boundary condition and ``do nothing" condition may be used. For a
flow in a vessel with in/outlet a mixture of Dirichlet condition
$v=0$ on the wall and pressure conditions on the inlet/outlet is
used. But for the flow in a vessel with in/outlet and a free surface
a mixture of Dirichlet condition, a Neumann condition $2\nu
\varepsilon(v)n-pn=0$ and pressure conditions may be used.

There are vast literatures for the Stokes and
 Navier-Stokes problems with mixed boundary conditions and several
 variational formulations are used for them, but two possible examples above are not considered except \cite{kc}.
With exception \cite{kc} mixtures of boundary conditions for
Navier-Stokes equations may be divided into three groups according
to what bilinear form is used for a variational formulation (for
more concrete one can refer to Introduction of \cite{kc}).

To include more different boundary conditions together, in \cite{kc}
the relations among strain, rotation, normal derivative of vector
field and shape of boundary surface are obtained and applied to the
stationary and non-stationary Navier-Stokes problems with mixture of
seven kinds of boundary conditions. However, for the non-stationary
Navier-Stokes problems we only were concerned with a mixed boundary
condition including total pressure (not static pressure), total
stress (not stress) and so on. Thus, in this paper we will study the
non-stationary problems with a mixed boundary condition including
static pressure (not total pressure), stress (not total stress) and
so on.

On the other hand, when one of static pressure (instead of total
pressure), stress (instead of total stress) or the outlet boundary
condition is given on a portion of boundary, for the initial
boundary value problems of the Navier-Stokes equations existence of
a unique local-in-time solution and a unique solution on a given
interval for small given data (in what follows we call it a solution
for small data) are proved. In the mathematical point of view the
main difficulty of such problems consists in the fact that in priori
estimation the inertial term is not canceled, and in the mechanical
point of view it is in the fact that the kinetic energy of fluid is
not controlled by the data and uncontrolled ``backward flow" can
take place at the portion of boundary(cf. preface in \cite{b}).

The Navier-Stokes equations with mixture of Dirichlet condition and
stress condition are studied. In \cite{ks} under smoothness
condition and a compatibility condition of data at the initial
instance existence of a unique local-in-time solution to the 3-D
Navier-Stokes equations is studied. In \cite{b} for the
Navier-Stokes equations on the polyhedral domain with mixture of
Dirichlet condition, Navier slip condition and stress condition a
local-in-time solution  is studied. Here smoothness of solutions to
the corresponding steady Stokes problem is used essentially.

The Navier-Stokes equations and the Boussinesq equations with
mixture of Dirichlet condition and the outlet boundary condition are
studied. For 2-D Navier-Stokes equations a local-in-time solution in
\cite{bk} and a solution for small data in \cite{bk1} are studied.
Here also smoothness of solutions to the corresponding 2-D steady
linear problem is important. For the Boussinesq systems a
local-in-time strong solutions in \cite{b0} on 2-D channel and in
\cite{b3} on 3-D channel are studied. Here smoothness of solutions
to the corresponding steady linear problems, respectively, in
\cite{b2} and \cite{b4} is the key.  In \cite{k} it is proved that
if under a compatibility condition at initial instance there exists
a unique solution, then so does for small perturbed data. This
result shows that under the compatibility condition there exists a
unique solution for small data.  In \cite{sk} for the Boussinesq
equations it is proved that under a compatibility condition there
exists a unique local-in-time solution, which is similar to the
result in
 \cite{ks}.
Smoothness of initial function in the compatibility condition of
\cite{k} is stronger than one in \cite{ks} and \cite{sk}.

The 2-D Navier-Stokes equations with mixture of Dirichlet condition
and pressure is studied. In \cite{m} existence of a unique solution
for small data is proved.

The Navier-Stokes equations with mixture of Dirichlet condition,
outlet condition and tangent stress condition is studied. In
\cite{b1} existence of a unique solution for small perturbation data
of the given solution is studied. Here also smoothness of solutions
to the corresponding steady linear problem is the key.

  In the present paper as a continuation of \cite{kc}, we are
concerned with the non-steady Navier-Stokes equations with mixed
boundary conditions involving  conditions for Dirichlet,  static
pressure, rotation, stress and normal derivative of velocity
together. Owing to the relations among strain, rotation, normal
derivative of velocity and shape of boundary surface obtained in
\cite{kc} (Theorems \ref{t1.1}, \ref{t1.2}), we can consider all
these boundary conditions together.

In general, the solution of the Stokes problem with mixed boundary
conditions has singularities on the intersections of surfaces for
different boundary conditions and the leading singular exponent of
the solution is
 a function of the intersection angle (cf. \cite{os}).
For the problem with Dirichlet condition and ``do nothing" condition
if the intersection angle is $\pi/2,$ then under some conditions for
data the solution belongs to $H^2(\Omega)$ (cf. \cite{b2}). For the
problem with Dirichlet condition and stress conditions, for similar
results refer to subsection 5.5 of \cite{mr} and section 10.3 of
\cite{mr1}.
 In our case the boundary conditions are more complicated than others,
 and there is no result for smoothness of solutions to the corresponding steady linear
problems. Thus, we prove
 existence of a unique weak solution for small data under a
compatibility condition at initial instance. We also prove that if a
solution smooth as in \cite{k} is given, then under the
compatibility condition for the small perturbed data there exists a
unique solution.

We are concerned with two problems distinguished according to
boundary conditions. Using relations among strain, rotation, normal
derivative of vector field and shape of boundary surface, which are
obtained in \cite{kc}, we reflect all these boundary conditions into
variational formulations for problems.

 This paper consists of 5 sections.
 In the end of this section the method in this paper is compared with
 another one.

In Section 2, the formulations of problems and some results for
definitions of weak solutions are given. According to bilinear forms
used for variational formulations for problems, the involved
boundary conditions are slightly different. Thus, difference between
our two problems is explained (Remark \ref{r1.0}).

 In Section 3, first, for the Navier-Stokes problems with seven kinds of boundary conditions a variational formulation, which is  based on the bilinear form
\begin{equation}\label{0.1}
 a(v,u)=2\Sigma_{i,j}(\varepsilon_{ij}(v),\varepsilon_{ij}(u))_{L_2(\Omega)}\quad
\mbox{for}\,\, v,u \in \textbf{H}^1(\Omega),
\end{equation}
is given. Next, by a transformation of the unknown function, the
problem is reduced to an equivalent problem in which the linear main
operator is positive definite. Then, studying properties of linear
operator differential equations and using a local diffeomorphism
theorem of nonlinear operator, we prove that under a compatibility
condition similar to one in \cite{k}, \cite{ks}, \cite{sk} there
exists a unique solution for small data (Theorem \ref{t2.1}).

In Section 4 for the Navier-Stokes problems with six kinds of
boundary conditions, which is a little different from one in Section
3, a variational formulation based on the bilinear form
\begin{equation}\label{0.2}
a(v,u)=(\nabla v,\nabla u)_{\textbf{L}_2(\Omega)}\quad
\mbox{for}\,\, v,u \in \textbf{H}^1(\Omega)
\end{equation}
 is given. Also,
by a transformation of the unknown function, the problem is reduced
to another equivalent problem in which the linear main operator is
positive definite. The result similar to one in Section 3 is
obtained (Theorem \ref{t3.1}).

Section 5 is considered in comparison with \cite{k} rather than
practical models. Existence of a unique solution for the small data
perturbed from a given solutions is proved under a compatibility
condition (Theorem \ref{t4.1}).

The compatibility conditions in Sections 3, 4 and 5 are similar to
one in  \cite{ks}, \cite{sk} and \cite{k}.  In point of view of
smoothness of the initial functions, the conditions are the same
with one in \cite{ks}, \cite{sk} concerning with local-in-time
solutions (cf. Remark \ref{r2.3}), but weaker than one in  \cite{k}
concerning with solutions for small data as our case (cf. Remarks
\ref{r3.2}, \ref{r4.1}). In \cite{k} the main results for the
nonlinear problem as perturbation of a linear problem is obtained by
a local diffeomorphism theorem relying on the properties of the
corresponding linear problem, and so is it in our paper.

Then, let us
consider why smoothness of the initial functions in our
compatibility conditions is weaker than one in \cite{k}.

 Let $\textbf{H}^k=(W^k_2(\Omega))^l $ be a Sobolev spaces on $\Omega$ with dimension $l,$ $V$ be a divergence-free subspace
 of $\textbf{H}^1$ satisfying appropriate boundary conditions, $H$ - the closure of $V$
 in $(L^2(\Omega))^l$,
 $V^*$- the adjoint space of $V$, $V^{r_0}(\Omega)=V\cap \textbf{H}^{r_0}(\Omega),$  where $r_0>l/2,$  $\mathscr{X}=\{w\in L_2(0,T;V); w'\in L_2(0,T;V),
w''\in L_2(0,T;V^*)\}, \mathscr{Y}=\{w\in L_2(0,T;V^*); w'\in
L_2(0,T;V^*)\}$ and $A:V\rightarrow V^*$- the Stokes operator.

Considering a linear problem
\[
 \left\{
 \begin{aligned}
    &u'(t)+Au(t)=f(t),\\
    &u(0)=\varphi,
 \end{aligned}
 \right.
\]
in \cite{k} the author proved the fact that a map
$u\rightarrow\{u(0),Lu\equiv u'+Au\}$ is linear continuous
one-to-one from ${\mathcal{X}}=\{u\in \mathscr{X}:u(0)\in
V^{r_0}(\Omega)\}$ onto ${\mathcal{Y}}=\{[\varphi,h]: \varphi\in
V^{r_0}(\Omega), h\in \mathscr{Y}, h(0)-A\varphi\in H\}$ (Theorem
3.1 in \cite{k}). Then, starting from this fact, the author studied
a nonlinear problem
\begin{equation}\label{0.3}
\left\{
 \begin{aligned}
    &u'(t)+(A+B)u(t)=f(t),\\
    &u(0)=\varphi,
 \end{aligned}
 \right.
\end{equation}
where $B: V\rightarrow V^*$ is defined by
$\left<Bu,v\right>=\left<(u\cdot\nabla)u,v\right>$ for $u,v\in V$.
 To this end, it was proved that
the inverse of a nonlinear map $u\rightarrow\{u(0),\tilde{L}u\equiv
u'+(A+B)u\}$ is one-to-one from a neighborhood of $0_{\mathcal{Y}}$
onto a neighborhood of $0_{\mathcal{X}}.$ From this fact the author
obtained that under the compatibility condition $f(0)-A\varphi\in H,
\varphi\in V^{r_0}(\Omega)$ and smallness of data, there exists a
unique solution to \eqref{0.3}.

However, we prove that for a modified operator $A$ a map
$u\rightarrow\{u'(0),Lu\equiv u'+Au\}$  is linear continuous
one-to-one from $\mathscr{X}$ onto $H\times\mathscr{Y}$ (Lemmas
\ref{2.3}, \ref{3.3}). Then using this fact, we prove that for a
modified operator $B(t)$ the inverse of a nonlinear map
$u\rightarrow\{u'(0),\tilde{L}u\equiv u'+(A+B(t))u\}$ is one-to-one
from a neighborhood of $0_{H\times\mathscr{Y}}$ onto a neighborhood
of $0_{\mathscr{X}}.$ By this  under the compatibility condition
$f(0)-A\varphi-B(0)\varphi \in H$ without $\varphi\in
V^{r_0}(\Omega)$ and smallness of data, we get existence of unique
solution to \eqref{0.3} (Theorems \ref{t2.1}, \ref{t3.1},
\ref{t4.1}, \ref{t4.2}). Since $B(0)\varphi=(\varphi\cdot
\nabla)\varphi,$ for $\varphi\in V^{l/2}(\Omega)$ we get
$B(0)\varphi\in H$, and so our condition is weaker than one in
\cite{k}.

 \section{Problems and preliminaries} \setcounter{equation}{0}

Throughout this paper we will use the following notation.

Let $\Omega$ be a connected bounded open subset of $ R^l,\,l=2,3.$
$\partial\Omega\in C^{0,1}$, $\partial\Omega=\cup_{i=1}^N
\overline{\Gamma}_i$,\,$\Gamma_i\cap\Gamma_j=\varnothing$ for $i\neq
j$, $\Gamma_i\in C^2$ for $i=2,3$. For Problems I and II stated
below we assume, respectively, $\Gamma_7\in C^2$ and $\Gamma_5\in
C^2$. Let $n(x)$ and $\tau(x)$ be, respectively, outward normal and
tangent unit vectors
 at $x$ in $\partial\Omega$. When $X$ is a Banach space,  $\mathbf{X}=X^l$ and  $\mathbf{X}^*$ is the dual of $\mathbf{X}$. Let $W^k_\alpha(\Omega)$ be
Sobolev spaces, $H^k(\Omega)=W^k_2(\Omega)$, and so
$\mathbf{H}^1(\Omega)=\{H^1(\Omega)\}^l$.  $Q=\Omega\times (0,T)$,
$\Sigma_i=\Gamma_i\times (0,T)$, $0<T<\infty$.

An inner product and a norm in the space $\mathbf{L}_2(\Omega)$ are
denoted, respectively, by $(\cdot\,,\cdot)$ and $\|\cdot\|$; and
$\langle\cdot\,,\cdot\rangle$ means the duality pairing between  a
Sobolev space $X$ and its dual one. Also,
$(\cdot\,,\cdot)_{\Gamma_i}$ is an inner product in the
$\mathbf{L}_2(\Gamma_i)$ or $L_2(\Gamma_i)$; and
$\langle\cdot\,,\cdot\rangle_{\Gamma_i}$ means the duality pairing
between $\mathbf{H}^{\frac{1}{2}}(\Gamma_i)$ and
$\mathbf{H}^{-\frac{1}{2}}(\Gamma_i)$ or between
$H^{\frac{1}{2}}(\Gamma_i)$ and $H^{-\frac{1}{2}}(\Gamma_i)$. The
inner product and norms in $R^l$, respectively, are denoted by
$(\cdot\, , \cdot)_{R^l}$ and $|\cdot|$. Sometimes $a\cdot b$ is
used for inner product in $R^l$ between $a$ and $b$. When
$\mathbf{X}$ is a Banach space, the zero element of $\mathbf{X}$ is
denoted by $0_\mathbf{X}$  and $\mathscr{O}_M(0_\mathbf{X})$ means
$M$-neighborhood of $0_\mathbf{X}$.

In this paper for the Navier-Stokes problem
\begin{equation}\label{1.1}
\left\{\begin{aligned}
 &\f{\p v}{\p t}-\nu\Delta v+(v\cdot\nabla)v+\nabla p=f,\,\,\,\quad \mbox{in}\,\,
 Q,\\&
 \nabla \cdot v=0,\\&v(0)=v_0
\end{aligned}\right.
\end{equation}
we are concerned with the problems I and II, which are distinguished
according to boundary conditions. Problem I is one with the boundary
conditions
\begin{equation}\label{1.2}
\begin{aligned}
&(1)\quad v|_{\Gamma_1}=h_1,\\
&(2)\quad v_\tau|_{\Gamma_2}=0,\,\,  -p|_{\Gamma_2}=\phi_2,\\
&(3)\quad v_n|_{\Gamma_3}=0,\,\,\mbox{rot}\,v\times n|_{\Gamma_3}=\phi_3/\nu,\\
&(4)\quad v_\tau|_{\Gamma_4}=h_4,\,\,(-p+2\nu\varepsilon_{nn}(v))|_{\Gamma_4}=\phi_4,\\
&(5)\quad v_n|_{\Gamma_5}=h_5,\,\,2(\nu\varepsilon_{n\tau}(v)+\alpha v_\tau)|_{\Gamma_5}=\phi_5,\\
&(6)\quad
(-pn+2\nu\varepsilon_n(v))|_{\Gamma_6}=\phi_6,\\
&(7)\quad v_\tau|_{\Gamma_7}=0,\,(-p+\nu\f{\p v}{\p n}\cdot
n)|_{\Gamma_7}=\phi_7,
\end{aligned}
\end{equation}
and Problem II is one with the conditions
\begin{equation}\label{1.3}
\begin{aligned}
&(1)\quad v|_{\Gamma_1}=h_1,\\
&(2)\quad v_\tau|_{\Gamma_2}=0,\,\,  -p|_{\Gamma_2}=\phi_2,\\
&(3)\quad v_n|_{\Gamma_3}=0,\,\,\mbox{rot}\,v\times n|_{\Gamma_3}=\phi_3/\nu,\\
&(4)\quad v_\tau|_{\Gamma_4}=h_4,\,\,(-p+2\nu\varepsilon_{nn}(v))|_{\Gamma_4}=\phi_4,\\
&(5)\quad v_n|_{\Gamma_5}=h_5,\,\,2(\nu\varepsilon_{n\tau}(v)+\alpha v_\tau)|_{\Gamma_5}=\phi_5,\\
&(6)\quad  (-pn+\nu\f{\p v}{\p n})|_{\Gamma_7}=\phi_7,
\end{aligned}
\end{equation}
together $\Gamma_6=\varnothing$. Here and in what follows
$u_n=u\cdot n$, $u_\tau=u-(u\cdot n)n$,
 $\varepsilon(v)=\{\varepsilon_{ij}(v)\},$  $\varepsilon_n(v)=\varepsilon(v)n$,
$\varepsilon_{nn}(v)=(\varepsilon(v)n,n)_{R^l}$,
$\varepsilon_{n\tau}(v)=\varepsilon(v)n-\varepsilon_{nn}(v)n$, and
 $h_i,\,\phi_i$ are functions or vector functions of $x,t$ defined on $\Gamma_i\times (0,T)$.

\begin{remark}\label{r1.0}
 The condition (6) of \eqref{1.3} (with $\phi_7=0$) is ``do nothing" condition, but
 (7) of \eqref{1.2} (with $\phi_7=0$) is rather
different from ``do nothing" condition and we can not unify two
problems.

First, let us consider why (7) of \eqref{1.2} is not changed with
(6) of \eqref{1.3}. In Section 3 relying on the bilinear form
\eqref{0.1} and integrating by parts $(-\nu\Delta v+\nabla p, u)$,
we get boundary integral $(-2\nu(\varepsilon(v)n, u)_{\p \Omega}+(p,
u\cdot n)_{\p \Omega}.$ Then,  in order to reflect the boundary
conditions into Formulation \ref{f2.1}, using $v_\tau=0$ or $v_n=0$
and applying Theorems \ref{t1.1} or \ref{t1.2}, we transform the
boundary integrals on $\Gamma_i, i=2,3,7.$ (cf.
\eqref{2.1}-\eqref{2.3}). Concretely, under conditions
$v_\tau|_{\Gamma_7}=0$
 we have
\begin{equation}\label{1.11}
\begin{aligned}
&\big<-pn+\nu\f{\p v}{\p n},u\big>_{\Gamma_7}\quad \forall  u
\,\,\mbox{with}\,\, u_\tau=0.
\end{aligned}
\end{equation}
 Usually, $v_\tau=0$ does not imply
$\f{\p v}{\p n}\cdot \tau=0,$ but
 by virtue of the conditions $u_\tau=0$ and (7) of \eqref{1.2} we have
\begin{equation}\label{1.13}
\begin{aligned}
&\big<-pn+\nu\f{\p v}{\p n},u\big>_{\Gamma_7}=\big<-p+\nu\f{\p v}{\p
n}n,u_n\big>_{\Gamma_7}=\left <\phi_7, u_n\right>_{\Gamma_7}\quad
\forall  u \,\,\mbox{with}\,\, u_\tau=0.
\end{aligned}
\end{equation}
Thus, substituting $\big<-pn+\nu\f{\p v}{\p n},u\big>_{\Gamma_7}$
with $\left <\phi_7, u_n\right>_{\Gamma_7}$, we reflect the boundary
condition (7) of \eqref{1.2} into Formulation \ref{f2.1}.

 Changing (7) of \eqref{1.2} by  $(-pn+\nu\f{\p v}{\p
n})|_{\Gamma_7}=\phi_7$ with a vector $\phi_7$ and substituting
$\big<-pn+\nu\f{\p v}{\p n},u\big>_{\Gamma_7}$ with $\left <\phi_7,
u\right>_{\Gamma_7}$, we can come to a formal variational
formulation. But when a solution $v$ is smooth enough, trying to
convert from the formal variational formulation to the original
problem, we come to
\begin{equation}\label{1.14}
(-pn+\nu\f{\p v}{\p n},
u)_{\Gamma_7}=\left<\phi_7,u\right>_{\Gamma_7}\quad \forall u,
u_\tau=0
\end{equation}
on $\Gamma_8.$  If we have \eqref{1.14} without $u_\tau=0$, then
from \eqref{1.14} we can get $-pn+\nu\f{\p v}{\p n}=\phi_7$ on
$\Gamma_8.$
 But owing to $u_\tau=0$ we only get $(-pn+\nu\f{\p v}{\p n}, n)_{\Gamma_7}=<\phi_7,n>_{\Gamma_7}$.
  This shows that the formal variational formulation is not equivalent to the original condition on $\Gamma_7$ and
  equivalent to $(-p+\nu\f{\p v}{\p n}n)|_{\Gamma_7}=\phi_7\cdot n.$ (Thus, (7) of (3.3) in \cite{kc} was
corrected. See Erratum to: \cite{kc}.)

Similarly, relying on the form \eqref{0.2}, we can reflect ``do
nothing" condition into Formulation \ref{f3.1}, but can not do (6)
of \eqref{1.2}. (cf. \eqref{3.1}-\eqref{3.3}). Therefore, two
problems are not unified.

  ``Do nothing" boundary condition  results from variational formulation based on \eqref{0.2} and
does not have a real physical meaning, but is rather used in
truncating large physical domains to smaller computational domains
by assuming parallel flow. If the flow is parallel in a near the
boundary, then (7) of \eqref{1.2} is same with ``do nothing"
condition.
\end{remark}\vspace*{.1cm}

For variational formulations of Problems I, II we need the
following.

 Let $\Gamma$ be a surface (curve for $l=2$) of $C^2$ and
$v$ be  a vector field of $C^2$ on a domain of $R^l$ near $\Gamma$.
In what follows the surface is a piece of boundary of 3-D or 2-D
bounded connected domains, and so we can assume the surface is
oriented.
\begin{theorem}\label{t1.1} (Theorem 2.1 in \cite{kc})
Suppose that $v\cdot n|_{\Gamma}=0$. Then, on the surface $\Gamma$
the followings hold.
\begin{equation}\label{1.4}
 \left(\varepsilon(v) n,
\tau\right)_{R^l}=\f{1}{2}(\mbox{rot}\,v\times n,
\tau)_{R^l}-(S\tilde{v},\tilde{\tau})_{R^{l-1}},
\end{equation}
\begin{equation}\label{1.5}
    (\mbox{rot}\,v\times n, \tau)_{R^l}=\left(\f{\p v}{\p n},
\tau\right)_{R^l}+(S\tilde{v},\tilde{\tau})_{R^{l-1}},
\end{equation}
\begin{equation}\label{1.6}
    \left(\varepsilon(v) n, \tau\right)_{R^l}=\f{1}{2}\left(\f{\p
v}{\p n}, \tau\right)_{R^l}-\f{1}{2}
(S\tilde{v},\tilde{\tau})_{R^{l-1}},
\end{equation}
where $S$ is the shape operator of the surface $\Gamma$ for $l=3$,
i.e.
$$
\begin{aligned}
&\hspace*{4cm} S=
\begin{pmatrix}
L & K\\
M & N
\end{pmatrix},\\&
L=\left(e_1, \f{\p n}{\p e_1}\right)_{R^l},\,\,K=\left(e_2, \f{\p
n}{\p e_1}\right)_{R^l},\,\, M=\left(e_1, \f{\p n}{\p
e_2}\right)_{R^l},\,\, N=\left(e_2, \f{\p n}{\p e_2}\right)_{R^l},
\end{aligned}
$$
and the curvature of $\Gamma$ for $l=2$. Here $e_i$ are the unit
vectors in a local curvilinear coordinates on $\Gamma$ and
$\tilde{v}, \tilde{\tau}$ are expressions of the vectors $v, \tau$
in the coordinate system.
\end{theorem}

\begin{theorem}\label{t1.2} (Theorem 2.2 in \cite{kc})
If $v_\tau|_{\Gamma}=0$ and $\text{div}\, v=0$, then on the surface
$\Gamma$ the following holds.
$$
\left(\varepsilon(v) n, n\right)_{R^l}=\left(\f{\p v}{\p n},
n\right)_{R^l}=-(k(x)v,n)_{R^l}
$$
where $k(x)=\mbox {div}\,n(x)$.
\end{theorem}

\begin{remark}\label{r1.1} (cf. \cite{kc}) $k(x)=\mbox {div}\,n(x)=Tr(S(x))=2\times \mbox{mean curvature}.$

If $\Gamma$ is a piece of  $\partial\Omega$, then since
$\partial\Omega\in C^{0,1}$ and $\Gamma\in C^2$, elements of $S$
belong to $C(\bar{\Gamma})$ and so does $k(x)$.
\end{remark}\vspace*{.1cm}

\section{Existence of a unique solution to problem I} \setcounter{equation}{0}

We use the following notation.

$\mathbf{V}=\{u\in \mathbf{H}^1(\Omega):\mbox{div}\,u=0,\,
u|_{\Gamma_1}=0,\, u_\tau|_{\Gamma_2\cup\Gamma_4\cup\Gamma_7}=0,\,
u_n|_{\Gamma_3\cup\Gamma_5}=0\}$ and
$\mathbf{V}_{\Gamma237}(\Omega)=\{u\in
\mathbf{H}^1(\Omega):\mbox{div}\,u=0,\,
u_\tau|_{\Gamma_2\cup\Gamma_7}=0,\, u_n|_{\Gamma_3}=0\}.$ Denote by
$H$ the completion of $\mathbf{V}$ in the space
$\textbf{L}_2(\Omega)$. Through this paper
$\widetilde{\textbf{V}}=\{u\in
\mathbf{H}^1(\Omega):\mbox{div}\,u=0\}$.

 By Theorems \ref{t1.1} and
\ref{t1.2} we have that for $v\in
\mathbf{H}^2(\Omega)\cap\mathbf{V}_{\Gamma237}(\Omega)$ and $u\in
\mathbf{V}$
\begin{equation}\label{2.1}
\begin{aligned}
 -(\Delta v, u)&=2(\varepsilon(v),
\varepsilon(u))-2(\varepsilon(v)n, u)_{\cup_{i=2}^7\Gamma_i}\\&
=2(\varepsilon(v),
\varepsilon(u))+2(k(x)v,u)_{\Gamma_2}-(\mbox{rot}\,v\times n,
u)_{\Gamma_3}+2(S\tilde{v},\tilde{u})_{\Gamma_3}\\&
\quad-2(\varepsilon_n(v),u)_{\cup_{i=4}^7\Gamma_i}\\&
=2(\varepsilon(v),
\varepsilon(u))+2(k(x)v,u)_{\Gamma_2}-(\mbox{rot}\,v\times u,
u)_{\Gamma_3}+2(S\tilde{v},\tilde{u})_{\Gamma_3}\\&
\quad-2(\varepsilon_{nn}(v),u\cdot
n)_{\Gamma_4}-2(\varepsilon_{n\tau}(v),u)_{\Gamma_5}-2(\varepsilon_n(v),u)_{\Gamma_6}\\&\quad-\left(\f{\p
v}{\p n},u\right)_{\Gamma_7}+(k(x)v,u)_{\Gamma_7}.
\end{aligned}
\end{equation}
Also, for $p\in H^1(\Omega)$ and $u\in \mathbf{V}$ we have
\begin{equation}\label{2.2}
\begin{aligned}
(\nabla p, u)&=(p, u\cdot n)_{\cup_{i=2}^7\Gamma_i}=(p, u\cdot
n)_{\Gamma_2}+(p, u\cdot n)_{\Gamma_4}+(pn,
u)_{\Gamma_6\cup\Gamma_7},
\end{aligned}
\end{equation}
where the fact that $u_n|_{\Gamma_3\cup\Gamma_5}=0$ was used.

Let
$$
\begin{aligned}
&\mathscr{X}=\{w\in L_2(0,T;\textbf{V}); w'\in L_2(0,T;\textbf{V}),
w''\in L_2(0,T;\textbf{V}^*)\},\\&
\hspace{1cm}\|w\|_{\mathscr{X}}=\|w\|_{L_2(0,T;\textbf{V})}+\|w'\|_{L_2(0,T;\textbf{V})}+\|w''\|_{L_2(0,T;\textbf{V}^*)},\\
&\mathscr{Y}=\{w\in L_2(0,T;\textbf{V}^*); w'\in
L_2(0,T;\textbf{V}^*)\},\\&
\hspace{1cm}\|w\|_{\mathscr{Y}}=\|w\|_{L_2(0,T;\textbf{V}^*)}+\|w'\|_{L_2(0,T;\textbf{V}^*)},\\&
\mathscr{W}=\{w\in L_2(0,T;\widetilde{\textbf{V}}); w'\in
L_2(0,T;\widetilde{\textbf{V}}), w''\in L_2(0,T;{\widetilde
{\textbf{V}}}^*)\},\\&
\hspace{1cm}\|w\|_{\mathscr{W}}=\|w\|_{L_2(0,T;\widetilde{\textbf{V}})}+\|w'\|_{L_2(0,T;\widetilde{\textbf{V}})}+
\|w''\|_{L_2(0,T;{\widetilde {\textbf{V}}}^*)}.
\end{aligned}
$$
Here and in what follows $w'$ means the derivative of $w(t)$ with
respect to $t$.

For Problem I, we use the following assumptions.
\begin{assumption}\label{a2.1}
 $f, \,f'\in L_2(0,T;\mathbf{V}^*),\,\, \phi_i,\,
\phi_i' \in L_2(0,T; H^{-\f{1}{2}}(\Gamma_i)),\,i=2,4,7,$ $\phi_i,
\,\phi_i'\in L_2(0,T;$ $
\mathbf{H}^{-\f{1}{2}}(\Gamma_i)),i=3,5,6,\,\,\alpha_{ij}\in
 L_\infty(\Gamma_5)$, where $\alpha_{ij}$  are components of the matrix $\alpha$, and
 $\Gamma_1\neq\emptyset$.
\end{assumption}\vspace*{.1cm}

\begin{assumption}\label{a2.2} There exists a function $U\in \mathscr{W}$ such that
$$
\begin{aligned}
\mbox{div}\,U=0,\,U|_{\Gamma_1}=h_1,\, U_\tau|_{\Gamma_2\cup
\Gamma_7}=0,\, U_n|_{\Gamma_3}=0,\, U_\tau|_{\Gamma_4}=h_4,\,
U_n|_{\Gamma_5}=h_5.
\end{aligned}
$$
Also, $U(0,x)-v_0\in \textbf{V}$.
\end{assumption}\vspace*{.1cm}

\begin{remark}\label{r3.0}  In practical situations $h_4, h_5=0,$
and in the cases if  for every fixed $t$ $h_1(t,x)\in
H^{\f{1}{2}}_{00}(\Gamma_1),$ $\int_{\Gamma_1} h_1(t, x)\cdot
n\,dx=0$ and  $\|h_1(t,x)\|_{ H^{\f{1}{2}}(\Gamma_1)}$ is smooth
enough with respect to $t$, then there exists such a function $U$.

\end{remark}\vspace*{.1cm}

Taking  \eqref{2.1} and \eqref{2.2} into account, we get the
following variational formulation for Problem I:
\begin{formulation}\label{f2.1}
Find $v$ such that
\begin{equation}\label{2.3}
\begin{aligned}
& v-U\in L_2(0,T;\mathbf{V}),\\& v(0)=v_0,\\& \langle v',u\rangle+
2\nu(\varepsilon(v),\varepsilon(u))+\langle(v\cdot\nabla)v,u\rangle+2\nu(k(x)v,u)_{\Gamma_2}\\&\hspace*{1cm}+2\nu(S\tilde{v},\tilde{u})_{\Gamma_3}+2(\alpha(x)v,u)_{\Gamma_5}+\nu(k(x)v,u)_{\Gamma_7}\\&
\hspace*{1cm}=\langle
f,u\rangle+\sum_{i=2,4,7}\langle\phi_i,u_n\rangle_{\Gamma_i}+\sum_{i=3,5,6}\langle\phi_i,u\rangle_{\Gamma_i}\quad
  \mbox{for all}\,\, u\in \textbf{V}.
\end{aligned}
\end{equation}
\end{formulation}

Taking Assumption \ref{a2.2} into account, put $v=\overline{z}+U$.
Then, we have the following problem equivalent to Formulation
\ref{f2.1}:

Find $\overline{z}$ such that
 \begin{equation}\label{2.4}
\begin{aligned}
& \overline{z}\in
L_2(0,T;\mathbf{V}),\\&\overline{z}(0)=\overline{z}_0\equiv
v_0-U(0)\in \mathbf{V},\\&
 \langle \overline{z}',u\rangle+
2\nu(\varepsilon(\overline{z}),\varepsilon(u))+\langle
(\overline{z}\cdot \nabla)\overline{z},u\rangle+\langle (U\cdot
\nabla)\overline{z},u\rangle+\langle (\overline{z}\cdot
\nabla)U,u\rangle\\&
\hspace{1cm}+2\nu(k(x)\overline{z},u)_{\Gamma_2}+2\nu(S\tilde{\overline{z}},\tilde{u})_{\Gamma_3}+2(\alpha(x)\overline{z},u)_{\Gamma_5}+\nu(k(x)\overline{z},u)_{\Gamma_7}\\&
 =-(U',u)-2\nu(\varepsilon(U),\varepsilon(u))-\langle (U\cdot
\nabla)U,u\rangle-2\nu(k(x)U,u)_{\Gamma_2}\\&\hspace{1cm}-2\nu(S\tilde{U},\tilde{u})_{\Gamma_3}
-2(\alpha(x)U,u)_{\Gamma_5}-\nu(k(x)U,u)_{\Gamma_7}+\langle
f,u\rangle\\&\hspace{1cm}+\sum_{i=2,4,7}\langle
\phi_i,u_n\rangle_{\Gamma_i}+\sum_{i=3,5,6}\langle\phi_i,u\rangle_{\Gamma_i}\qquad\mbox{for
all}\,\, u\in \mathbf{V}.
\end{aligned}
\end{equation}

Now, define an operator $A_0:\textbf{V}\rightarrow \textbf{V}^*$ by
 \begin{equation}\label{2.5}
\begin{aligned}
 \langle A_0y,u\rangle=&
2\nu(\varepsilon(y),\varepsilon(u))+2\nu(k(x)y,u)_{\Gamma_2}+2\nu(S\tilde{y},\tilde{u})_{\Gamma_3}\\&+2(\alpha(x)y,u)_{\Gamma_5}+\nu(k(x)y,u)_{\Gamma_7}\quad
\mbox{for all}\,\, y,u\in \mathbf{V}.
\end{aligned}
\end{equation}

\begin{lemma}\label{l2.1} $\exists\delta>0$,  $\exists k_0\geq 0$;
$ \langle A_0u,u\rangle\geq
\delta\|u\|_\textbf{\textbf{V}}^2-k_0\|u\|_H^2\quad\mbox{for
all}\,\, u\in \mathbf{V}. $
\end{lemma}
\emph{Proof}\,\, By Korn's inequality
\begin{equation}\label{2.6}
2\nu(\varepsilon(u),\varepsilon(u))\geq \beta
\|u\|^2_\mathbf{V}\quad\exists\beta>0, \mbox{for all}\,\,
 u\in \mathbf{V}.
\end{equation}
By Remark \ref{r1.1} and Assumption \ref{a2.1}, there exists a
constant $M$ such that
\[
\begin{aligned}
\|S(x)\|_\infty,\, \|k(x)\|_\infty, \|\alpha(x)\|_\infty \leq M,
\end{aligned}
\]
and so  there exists a constant $c_0$ (depending on $\beta$) such
that
\begin{equation}\label{2.7}
\begin{aligned}
&\left|2\nu(k(x)z,z)_{\Gamma_2}
+2\nu(S\tilde{z},\tilde{z})_{\Gamma_3}+\nu(k(x)z,z)_{\Gamma_7}+2(\alpha(x)y,u)_{\Gamma_5}\right|\\&
\hspace{2cm}\leq \f{\beta}{2}\|
z\|^2_{\mathbf{H}^1(\Omega)}+c_0\|z\|_{H}^2\,dt\quad \mbox{for
all}\,\, z\in \mathbf{V}
\end{aligned}
\end{equation}
((cf. Theorem 1.6.6 in \cite{bs} or (1), p. 258 in \cite{e})). Put
$\delta=\f{\beta}{2}$, $k_0=c_0$. Then, by \eqref{2.6}, \eqref{2.7}
we come to the asserted conclusion. $\square$\vspace*{.1cm}

\begin{remark}\label{r2.1}
In process of proof of Lemma \ref{l2.1}, we see that if
$\Gamma_i=\emptyset,\, i=2,3,7,$ or these are unions of pieces of
planes (segments in case of 2-D)and $\Gamma_5=\emptyset$ or
$\alpha(x)=0 $, then we can take $k_0=0$.
\end{remark}\vspace*{.1cm}

When $k_0>0$, if $k_0$ is not small enough, then the operator
defined by \eqref{2.5} is not positive, and so let us transform the
unknown function to get a positive operator $A$ in \eqref{2.9}
bellow. Now, let $k_0$ be the constant in Lemma \ref{l2.1} and put
$z=e^{-k_0t}\overline{z}$. Then, since
$e^{-k_0t}\overline{z}'=z'+k_0z$, we get the following problem
equivalent to problem \eqref{2.4}:

Find $z$ such that
 \begin{equation}\label{2.8}
\begin{aligned}
& z\in L_2(0,T;\mathbf{V}),\\&z(0)=v_0-U(0)\in \mathbf{V},\\&
 \langle z'(t),u\rangle+
2\nu(\varepsilon(z(t)),\varepsilon(u))+e^{k_0t}\langle (z(t)\cdot
\nabla)z(t),u\rangle+\langle (U(t)\cdot
\nabla)z(t),u\rangle\\&\hspace{1cm}+\langle (z(t)\cdot
\nabla)U(t),u\rangle+k_0(z(t),u)+2\nu(k(x)z(t),u)_{\Gamma_2}\\&\hspace{1cm}+2\nu(S\tilde{z}(t),\tilde{u})_{\Gamma_3}+2(\alpha(x)z(t),u)_{\Gamma_5}+\nu(k(x)z(t),u)_{\Gamma_7}\\&
 =e^{-k_0t}\Big[-(U'(t),u)-2\nu(\varepsilon(U(t)),\varepsilon(u))-\langle (U(t)\cdot
\nabla)U(t),u\rangle\\&\hspace{1cm}-2\nu(k(x)U(t),u)_{\Gamma_2}-2\nu(S\tilde{U}(t),\tilde{u})_{\Gamma_3}
-2(\alpha(x)U(t),u)_{\Gamma_5}\\&\hspace{1cm}-\nu(k(x)U(t),u)_{\Gamma_7}+\langle
f(t),u\rangle +\sum_{i=2,4,7}\langle
\phi_i(t),u_n\rangle_{\Gamma_i}\\&\hspace{1cm}+\sum_{i=3,5,6}\langle\phi_i(t),u\rangle_{\Gamma_i}\Big]\quad\mbox{for
all}\,\, u\in \mathbf{V}.
\end{aligned}
\end{equation}

Define operators $A$, $A_U(t):\textbf{V}\rightarrow \textbf{V}^*$ by
\begin{equation}\label{2.9}
\langle Av,u\rangle= \langle A_0v,u\rangle+(k_0v,u)\quad\mbox{for
all}\,\,
 v, u\in \mathbf{V},
\end{equation}
\begin{equation}\label{2.10}
\begin{aligned}
&\langle A_U(t)v,u\rangle=  \langle(U(t,x)\cdot
\nabla)v,u\rangle+\langle(v\cdot\nabla)U(t,x),u\rangle\quad\mbox{for
all}\,\,
 v, u\in \mathbf{V},
\end{aligned}
\end{equation}
where $A_0$ is the operator by \eqref{2.5} and $k_0$ is one in Lemma
\ref{l2.1}.  Since $U\in \mathscr{W}$, we have $U\in
C\left([0,T];\textbf{H}^1(\Omega)\right)$ and so such a definition
is well. Then, the operator $A$ is positive definite, and this fact
is used in future.

Define an operator $B(t): \textbf{V}\rightarrow \textbf{V}^*$ and
$F(t)\in V^*$ by
\begin{equation}\label{2.9.2} \langle
B(t)v,u\rangle=e^{k_0t}\langle (v\cdot\nabla)
v,u\rangle\quad\mbox{for all}\,\, v,u\in \mathbf{V},
\end{equation}
\begin{equation}\label{2.11}
\begin{aligned}
\langle
F(t),u\rangle&=e^{-k_0t}\Big[-(U'(t),u)-2\nu(\varepsilon(U)(t),\varepsilon(u))-\langle
(U(t)\cdot
\nabla)U(t),u\rangle\\&-2\nu(k(x)U(t),u)_{\Gamma_2}-2\nu(S\tilde{U}(t),\tilde{u})_{\Gamma_3}-2(\alpha(x)U(t),u)_{\Gamma_5}\\&-\nu(k(x)U(t),u)_{\Gamma_7}+\langle
f(t),u\rangle+\sum_{i=2,4,7}\langle
\phi_i(t),u_n\rangle_{\Gamma_i}\\&+\sum_{i=3,5,6}\langle\phi_i(t),u\rangle_{\Gamma_i}\Big]
\quad \mbox{for all}\,\, u\in \mathbf{V}.
\end{aligned}\end{equation}

Then, \eqref{2.8} is written by
 \begin{equation}\label{2.8.1}
\begin{aligned}
& z\in L_2(0,T;\mathbf{V}),\\&z(0)=v_0-U(0)\in \mathbf{V},\\&
 z'(t)+\left(A+A_U(t)+B(t)\right)z(t)=F(t).
\end{aligned}
\end{equation}

Now, define operators $L, {\widetilde A}_U, L_U, {\widetilde B}:
\mathscr{X}\rightarrow\mathscr{Y}$,
$C:\mathscr{X}\times\mathscr{X}\rightarrow\mathscr{Y}$ and $F\in
\mathscr{Y}$  by
 \begin{equation}\label{2.12.1}
\begin{aligned}
&\langle (Lz)(t),u\rangle=\langle z'(t),u\rangle+\langle
Az(t),u\rangle \quad\mbox{for all}\,\,
z\in \mathscr{X}, \mbox{for all}\,\, u\in \mathbf{V},\\
&\langle ({\widetilde A}_Uz)(t),u\rangle=\langle A_U(t)z(t),u\rangle
\quad\mbox{for all}\,\,
z\in \mathscr{X}, \mbox{for all}\,\, u\in \mathbf{V},\\
&\langle (L_Uz)(t),u\rangle=\langle z'(t),u\rangle+\langle
(A+A_U(t))z(t),u\rangle \quad\mbox{for all}\,\,
z\in \mathscr{X}, \mbox{for all}\,\, u\in \mathbf{V},\\
&\langle ({\widetilde B}z)(t),u\rangle=\langle
B(t)z(t),u\rangle\quad\mbox{for all}\,\,
z\in \mathscr{X}, \mbox{for all}\,\, u\in \mathbf{V},\\
&\langle C(w,z)(t),u\rangle=e^{k_0t}\langle (w(t)\cdot\nabla)
z(t),u\rangle+e^{k_0t}\langle
(z(t)\cdot\nabla)w(t),u\rangle\\&\hspace{3cm}\quad\mbox{for all}\,\,
w,z\in \mathscr{X}, \mbox{for all}\,\, u\in \mathbf{V},\\
 &(F)(t)= F(t).
\end{aligned}
\end{equation}

 \begin{lemma}\label{l2.2}  $C$ is a bilinear continuous operator such that $\mathscr{X}\times \mathscr{X}\rightarrow\mathscr{Y}$. Under Assumptions \ref{a2.1} and
 \ref{a2.2}, ${\widetilde A}_U$ is a linear continuous operator such that $\mathscr{X}\rightarrow\mathscr{Y}$ and $F\in \mathscr{Y}$.
\end{lemma}
\emph{Proof}\,\,Obviously, $C$ is bilinear. When $w\in \mathscr{X}$,
$$
w\in L_\infty(0, T; \textbf{V}),\quad\|w\|_{L_\infty(0, T;
\textbf{V})}\leq c\big [\|w\|_{L_2(0, T; \textbf{V})}+\|w'\|_{L_2(0,
T; \textbf{V})}\big]
$$
and by virtue of H\"{o}lder inequality and the imbedding theorem
$$
\begin{aligned}
&\big|e^{k_0t}\langle (w\cdot\nabla) z,u\rangle +e^{k_0t}\langle
(z\cdot\nabla)w,u\rangle\big|\\
&\hspace*{1cm}\leq c(\|w\|_{\textbf{L}_3} \|\nabla
z\|_{\textbf{L}_2}\|u\|_{\textbf{L}_6}+\|z\|_{\textbf{L}_3} \|\nabla
w\|_{\textbf{L}_2}\|u\|_{\textbf{L}_6})\leq c\|w\|_\textbf{V} \|
z\|_\textbf{V}\|u\|_\textbf{V}\\& \hspace*{4cm}\mbox{for all}\,\,
w,z,u\in
 \mathbf{V}.
\end{aligned}
$$
Thus,\begin{equation}\label{2.13}
\|C(w,z)\|_{L_2(0,T;\textbf{V}^*)}\leq
c\|w\|_{L_\infty(0,T;\textbf{V})} \|z\|_{L_2(0,T;\textbf{V})}\leq
c\|w\|_{\mathscr{X}}\cdot \|z\|_{\mathscr{X}}.
\end{equation}
Also, since
$$
\begin{aligned}
|\langle C(w,z)'(t),u\rangle|&=e^{k_0t}\big|k_0\langle
(w\cdot\nabla)z,u\rangle+k_0\langle (z\cdot\nabla)w,u\rangle+\langle
(w'\cdot\nabla)z,u\rangle\\&+\langle
(w\cdot\nabla)z',u\rangle+\langle (z'\cdot\nabla)w,u\rangle+\langle
(z\cdot\nabla)w',u\rangle\big|,
\end{aligned}
$$
taking \eqref{2.13} into account we have
\begin{equation}\label{2.14}
\begin{aligned}
\|C(w,z)'\|&_{L_2(0,T;\textbf{V}^*)}\leq
c\|C(w,z)\|_{L_2(0,T;\textbf{V}^*)}\\&+c\big[\big(\|w'\|_{L_2(0,T;\textbf{V})}+\|w\|_{L_\infty(0,T;\textbf{V})})(\|z'\|_{L_2(0,T;\textbf{V})}+\|z\|_{L_\infty(0,T;\textbf{V})}\big)\big]\\&
\leq c\|w\|_{\mathscr{X}}\cdot \|z\|_{\mathscr{X}}.
\end{aligned}
\end{equation}
\eqref{2.13} and \eqref{2.14} imply
\begin{equation}\label{2.15}
\begin{aligned}
\|C(w,z)\|_{\mathscr{Y}}\leq  c\|w\|_{\mathscr{X}}\cdot
\|z\|_{\mathscr{X}}.
\end{aligned}
\end{equation}

 By the same argument above, we have
\begin{equation}\label{2.16}
\begin{aligned}
\|{\widetilde A}_Uz\|_{\mathscr{Y}}\leq  c\|U\|_{\mathscr{W}}\cdot
\|z\|_{\mathscr{X}}.
\end{aligned}
\end{equation}

By Assumption \ref{a2.1}, Remark \ref{r1.1} and the trace theorem,
we can see that $F\in \mathscr{Y}$. $\square$\vspace*{.1cm}

\begin{lemma}\label{l2.3} The operator $\overline{L}$ defined by
$\overline{L}z=(z'(0), Lz)$ for $z\in \mathscr{X}$ is  a linear
continuous one-to-one operator from $\mathscr{X}$ onto $H\times
\mathscr{Y}$.
\end{lemma}
\emph{Proof}\,\, The linearity of $\overline{L}$ is obvious. The
fact $z\in \mathscr{X}$ implies that $z'\in C([0,T]; H)$,
$\|z'\|_{C([0,T]; H)}\leq c\|z\|_{\mathscr{X}}$, and so we see that
a map $z\in \mathscr{X}\rightarrow z'(0)\in H$ is continuous.

Clearly, $\|z'\|_{\mathscr{Y}}\leq c\|z\|_{\mathscr{X}}$. Also, by
Assumption \ref{a2.1}, Remark \ref{r1.1} and the trace theorem,
\begin{equation}\label{2.17}
|\langle Av,u\rangle|\leq c\|v\|_\mathbf{V}\cdot \|u\|_\mathbf{V}
\quad\mbox{for all}\,\,
 v,u\in \mathbf{V}.
\end{equation}
Formula \eqref{2.17} implies that the mapping
$z\in\mathscr{X}\rightarrow Az\in\mathscr{Y}$ is continuous.
Therefore, $\overline{L}$ is continuous.

Let us show that $\overline{L}$ is a one-to-one and surjective
operator from $\mathscr{X}$ onto $H\times \mathscr{Y}$.

First, let us prove that this operator is injective. For this, it is
enough to prove that the inverse image of $(0_H,0_\mathscr{Y})\in
H\times \mathscr{Y}$ by the operator $\overline{L}$ is
$0_\mathscr{X}$.
 By Lemma \ref{l2.1} and \eqref{2.9}, we get
\begin{equation}\label{2.18}
\langle Av,v\rangle\geq \delta
\|v\|^2_\mathbf{V}\quad\exists\delta>0, \mbox{for all}\,\,
 v\in \mathbf{V}.
\end{equation}
By \eqref{2.17}, \eqref{2.18} for any $q\in \textbf{V}^*$ there
exists a unique solution $y\in V$ to the following problem
\begin{equation}\label{2.19}
Ay=q.
\end{equation}
Let $z\in\mathscr{X}$ be the inverse image of
$(0_H,0_\mathscr{Y})\in H\times \mathscr{Y}$ by $\overline{L}$. Then
since $z'(0)=0_H$, putting $t=0$ in the first equation of
\eqref{2.12.1} we get
$$
\langle Az(0),u\rangle=0 \quad\mbox{for all}\,\,
 u\in \mathbf{V},
$$
where $z(0)=z(0,x)$. This means that $z(0)$ is a unique solution to
\eqref{2.19} for $q=0_{\textbf{V}^*}$, i.e. $z(0)=0_\textbf{V}$.
Putting $w=z'$, we get  $w(0)=z'(0)=0_H$. Taking $Lz=0$ into account
and differentiating the first equation of \eqref{2.12.1}, we have
\begin{equation}\label{2.20}
\begin{aligned}
&\langle w'(t),u\rangle+\langle Aw(t),u\rangle=0\quad\mbox{for
all}\,\,
 u\in \mathbf{V}.
 \end{aligned}
\end{equation}
The operator $A$ in \eqref{2.20} satisfies all conditions of Theorem
1.1, ch. 6 in \cite{ggz}. Thus, for problem \eqref{2.20} with an
initial condition $w(0)\in H$ there exists a unique solution $w$
such that $w\in L_2(0,T;\textbf{V}) ,\,\,w'\in
L_2(0,T;\textbf{V}^*)$. Since $w(0)=0_H$, we have $w=0$, which means
$z=0_{\mathscr{X}}$ since $z(0)=0_\textbf{V}$.

Let us prove that $\overline{L}$ is surjective. Let $(w_0,g)\in
H\times \mathscr{Y}$. Since $g\in \mathscr{Y}$, we have $g(0)\in
\textbf{V}^*$. Then, by \eqref{2.17} and \eqref{2.18}, there exists
a unique solution $z_0\in \textbf{V}$ to problem
\begin{equation}\label{2.21}
Az_0=g(0)-w_0.
\end{equation}
Let us consider problem
\begin{equation}\label{2.22}
\left\{\begin{aligned} &w'+Aw=g',\\& w(0)=w_0.
 \end{aligned}
\right.\end{equation} There exists a unique solution $w$ such that
$w\in L_2(0,T;\textbf{V}) ,\,\,w'\in L_2(0,T;\textbf{V}^*)$ to
problem \eqref{2.22} (cf. Theorem 1.3 of ch. 6 in \cite{ggz}). Now,
put
\begin{equation}\label{2.23}
 z=z_0+\int_0^tw(s)\,ds,
\end{equation}
where $ z_0$ is the solution to \eqref{2.21}. Then, $z'=w$ and $z\in
\mathscr{X}$. Integrating two sides of the first one in \eqref{2.22}
from $0$ to $t$ and using \eqref{2.23}, we have
\begin{equation}\label{2.24}
\begin{aligned}
&\langle w(t),u\rangle+\left\langle Az(t),u\right\rangle
-\left[\langle w_0,u\rangle+\langle Az_0,u\rangle \right]=\langle
g(t),u\rangle-\langle g(0),u\rangle\quad\mbox{for all}\,\,
 u\in \mathbf{V}.
\end{aligned}
\end{equation}
Taking  \eqref{2.21},  \eqref{2.23} into account, from \eqref{2.24}
we get
\begin{equation}\label{2.25}
\langle z'(t),u\rangle+\left\langle Az(t),u\right\rangle =\langle
g(t),u\rangle\quad\mbox{for all}\,\,
 u\in \mathbf{V}.
\end{equation}
This means that $z\in \mathscr{X}$ defined by \eqref{2.23} is the
inverse image of $(w_0,g)\in H\times \mathscr{Y}$ by the operator
$\overline{L}$, i.e. $\overline{L}$ is surjective. Therefore,
$\overline{L}$ is an epimorphism.
 $\square$\vspace*{.1cm}

\begin{lemma}\label{l2.4} Under Assumption \ref{a2.2}, let $\|U(0,x)\|_{\widetilde {\textbf{V}}}$ be small enough.
 The operator $\overline{L}_U$
defined by $\overline{L}_Uz=(z'(0), L_Uz)$ for $z\in \mathscr{X}$ is
a linear continuous one-to-one operator from $\mathscr{X}$ onto
$H\times \mathscr{Y}$.
\end{lemma}
\emph{Proof}\,\,When $z\in \mathscr{X}$, $z\in C([0, T];
\textbf{V})$ and
$$
\|z\|_{C([0, T]; \textbf{V})}\leq c\big [\|z\|_{L_2(0, T;
\textbf{V})}+\|z'\|_{L_2(0, T; \textbf{V})}\big].
$$ By virtue of this fact and Lemma \ref{l2.2}, the operator $\overline{A}_U\in (\mathscr{X}\rightarrow
H\times \mathscr{Y})$ defined by $\overline{A}_Uz=(0_H,{\widetilde
A}_Uz)$ is continuous. Thus, the operator $\overline{L}_U$ defined
on $\mathscr{X}$ is linear continuous.

As in Lemma 3.5 of \cite{k} it is proved that the operator
${\widetilde A}_U\in (\mathscr{X}\rightarrow\mathscr{Y})$ is
compact. Thus, $\overline{A}_U\in (\mathscr{X}\rightarrow
H\times\mathscr{Y})$ is also compact. Since
$\overline{L}_U=\overline{L}+\overline{A}_U$, by virtue of Theorem
3.4 in \cite{k} and Lemma \ref{l2.3} we know that in order to prove
that the operator $\overline{L}_U$ is one-to-one from $\mathscr{X}$
onto $H\times \mathscr{Y}$ it is enough to prove that
$\overline{L}_U$ is one-to-one from $\mathscr{X}$ into $H\times
\mathscr{Y}$. To prove the last fact it is enough to show that the
inverse image of $(0_H, 0_{\mathscr{Y}})$ by $\overline{L}_U$ is
$0_{\mathscr{X}}$. By H\"older inequality and imbedding theorem
\begin{equation}\label{2.27}
\begin{aligned}
\big| \big\langle(U(t,x)\cdot
\nabla)v,v\big\rangle+\big\langle(v\cdot\nabla)U(t,x),v\big\rangle\big|\leq
K_0 \|v\|_\textbf{V}\|U(t,x)\|_{\textbf{H}^1}\|v\|_\textbf{V}.
 \end{aligned}
\end{equation}
Thus, if $\|U(0,x)\|_{\widetilde V}$ is so small that
$\|U(0,x)\|_{\textbf{H}^1}\leq \f{\delta}{2K_0}$, then \eqref{2.17},
\eqref{2.18} and \eqref{2.27} imply
\begin{equation}\label{2.28}
 \big|\big\langle \big(A+A_U(0)\big)v,u\big\rangle\big|\leq c\|v\|_\mathbf{V}\cdot \|u\|_\mathbf{V}, \quad  \big\langle \big(A+A_U(0)\big)v,v\big\rangle\geq
\f{\delta}{2}\|v\|^2_{\mathbf{V}}\quad\mbox{for all}\,\,
 v,u\in \mathbf{V}.
\end{equation}
By \eqref{2.28} for any $q\in \textbf{V}^*$ there exists a unique
solution  $y\in V$ to
\begin{equation}\label{2.29}
(A+A_U(0))y=q.
\end{equation}
Let $z\in\mathscr{X}$ be the inverse image of $(0_H,0_\mathscr{Y})$
by $\overline{L}$. Then $z'(0)=0_H$, and putting $t=0$ from the
third one in \eqref{2.12.1} we get
$$
\big\langle \big(A+A_U(0)\big)z(0),u\big\rangle=0 \quad\mbox{for
all}\,\,
 u\in \mathbf{V},
$$
where $z(0)=z(0,x)$. This means that $z(0)$ is the unique solution
to \eqref{2.29} with $q=0_{\textbf{V}^*}$, i.e. $z(0)=0_\textbf{V}$.
Therefore, $z\in \mathscr{X}$ satisfies
\begin{equation}\label{2.30}
\left\{\begin{aligned}
&z'(t)+\big(A+A_U(t)\big)z(t)=0,\\
 &z(0)=0_\textbf{V}.
 \end{aligned}\right.
\end{equation}
Now, making duality pairing with $z(t)$ on two sides of $$
z'(t)+Az(t)=-A_U(t)z(t)
$$
and taking \eqref{2.18} into account and using Gronwall's
inequality, we can prove $z=0_{\mathscr{X}}$ as in Lemma 3.8 of
\cite{k}. It is finished to prove the Lemma. $\square$\vspace*{.1cm}

\begin{lemma}\label{l2.5} Under Assumption
 \ref{a2.2} the operator $T$ defined by $Tz=\big(z'(0), (L_U$ $+{\widetilde
B})z\big)$ for $z\in \mathscr{X}$ is continuously differentiable,
$T(0_\mathscr{X})=\left(0_H,0_\mathscr{Y}\right)$ and  the Frechet
derivative of $T$ at $0_\mathscr{X}$ is $\overline{L}_U$.
\end{lemma}
\emph{Proof}\,\, It is easy to verify that
$T(0_\mathscr{X})=\left(0_H,0_\mathscr{Y}\right)$. Since the
operator $L_U$ is linear, its Frechet derivative is the same with
itself. Therefore, if ${\widetilde B}$ is continuously
differentiable, then so is $T$.

For any $w, z\in \mathscr{X}$,
$$
\big({\widetilde B}(w+z)-{\widetilde B}w\big)(t)=
e^{k_0t}\big(w(t)\cdot\nabla\big)
z(t)+e^{k_0t}\big(z(t)\cdot\nabla\big)w(t)+({\widetilde B}z)(t).
$$
By \eqref{2.15}, we get
$$
 \lim_{\|z\|_{\mathscr{X}}\rightarrow
0}\f{\|{\widetilde B}z\|_{\mathscr{Y}}}{\|z\|_{\mathscr{X}}}\leq
\lim_{\|z\|_{\mathscr{X}}\rightarrow
0}\f{c\|z\|^2_{\mathscr{X}}}{\|z\|_{\mathscr{X}}}=0.
$$
Then, put
$$
C(w,z)(t)\equiv e^{k_0t}\big(w(t)\cdot\nabla\big)
z(t)+e^{k_0t}\big(z(t)\cdot\nabla\big)w(t)=({\widetilde B}'_w z)(t).
$$
By Lemma \ref{l2.2} ${\widetilde B}'_w \in (\mathscr{X}\rightarrow
\mathscr{Y})$ is continuous, and it is the Frechet derivative of
$\widetilde B$ at $w$ and also continuous with respect to $w$. Thus,
$T$ is continuously differentiable. Also from the formula above we
can see that the Frechet derivative of $\widetilde B$ at
$0_\mathscr{X}$ is zero. Therefore, the Frechet derivative of $T$ at
$0_\mathscr{X}$ is $\overline{L}_U$. $\square$\vspace*{.1cm}

Let us consider problem
\begin{equation}\label{2.31}
\big(A+A_U(0)+B(0)\big)u=q.
\end{equation}
\begin{lemma}\label{l2.6} Assume that $\|U(0,x)\|_{\widetilde {\textbf{V}}}$ is small enough.
If the norm of $q\in V^*$ is small enough, then there exists a
unique solution to \eqref{2.31} in some
$\mathscr{O}_M(0_\mathbf{V})$.
\end{lemma}
\emph{Proof}\,\, Since $\|U(0,x)\|_{\widetilde {\textbf{V}}}$ is
small enough, by \eqref{2.28},  for any fixed $z\in \mathbf{V}$
there exists a unique solution to problem
\begin{equation}\label{2.32}
\big(A+A_U(0)\big)w=q-B(0)z.
 \end{equation}
On the other hand,
\begin{equation}\label{2.33}
\begin{aligned}
|\langle B(0)w_1-B(0)w_2,u\rangle|\leq &K
M\|w_1-w_2\|_\mathbf{V}\cdot\|u\|_\mathbf{V}\quad\mbox{for all}\,\,
 w_i\in \mathscr{O}_M(0_\mathbf{V}),\,\, \mbox{for all}\,\, u\in
 \mathbf{V}.
 \end{aligned}
\end{equation}
Owing to \eqref{2.28} the solution  $w$ to \eqref{2.32} is estimated
as follows
\[
\|w\|_\mathbf{V}\leq
\f{2}{\delta}\big(\|q\|_{\mathbf{V}^*}+\|B(0)z\|_{\mathbf{V}^*}\big)\leq
\f{2}{\delta}\big(\|q\|_{\mathbf{V}^*}+KM^2\big).
\]
Thus, if $\|q\|_{\mathbf{V}^*}$ and $M$ are small enough, then
 the operator $(z\rightarrow w)$ maps $\mathscr{O}_M(0_\mathbf{V})$ into itself and
by \eqref{2.33} this operator is strictly contract. Therefore, in
$\mathscr{O}_M(0_\mathbf{V})$ there exists a unique solution to
\eqref{2.32}. Thus, we come to the asserted conclusion.
$\square$\vspace*{.1cm}

 For proof of unique existence of a solution to Problem I, we use the following
\begin{proposition}\label{p2.1} (cf. Theorem 4.1.1 in \cite{dm}) Let
$X, Y$ be Banach spaces, $\cal G$ an open set in $X$,
$f:X\rightarrow Y$ continuously differentiable on $\cal G$. Let the
derivative $f'(a)$ be an isomorphism of $X$ onto $Y$ for $a\in \cal
G$. Then there exist neighborhoods $\cal U$ of $a$, $\cal V$ of
$f(a)$such that $f$ is injective on $\cal U$, $f(\cal U)=\cal V$.
\end{proposition}

One of main results of this paper is the following
\begin{theorem}\label{t2.1} Suppose that Assumptions \ref{a2.1} and \ref{a2.2} hold.
 Assume that $\|U\|_{\mathscr{W}}$ and the norms of
 $f, f', \phi_i, \phi_i'$ in the spaces where they belong to are small enough.

If
\begin{equation}\label{2.34}
\begin{aligned}
w_0\equiv F(0)-(A+A_U(0)+B(0))z_0\in H,
\end{aligned}
\end{equation}
where $z_0=v_0-U(0,\cdot),$ and $\|w_0\|_H$ is small enough, then
there exists a unique solution to \eqref{2.3} in the space
$\mathscr{W}.$
\end{theorem}

\emph{Proof}\,\,First, let us prove existence of a solution.\\
If $\|U\|_{\mathscr{W}}$ and the norms of
 $f, f', \phi_i, \phi_i'$ in the spaces they belong to are small enough,
then $\|F\|_{\mathscr{Y}}$ is also small enough. By virtue of Lemmas
\ref{l2.4}, \ref{l2.5} and Proposition \ref{p2.1}, for any $R_1>0$
small enough if $\|F\|_{\mathscr{Y}}$, $R$ are small enough and
$w_1\in \mathscr{O}_R(0_H)$, there exists a unique $z\in
\mathscr{O}_{R_1}(0_{\mathscr{X}})$  such that
\begin{equation}\label{2.35}
\begin{aligned}
 z'(t)+\big(A+A_U(t)+B(t)\big)z(t)=F(t),\quad
z'(0)=w_1\in \mathscr{O}_R(0_H).
\end{aligned}
\end{equation}
Putting $t=0$ in \eqref{2.35}, we get
$$
F(0)-\big(A+A_U(0)+B(0)\big)z(0)=w_1\in\mathscr{O}_R(0_H).
$$

On the other hand, if $\|U\|_{\mathscr{W}}$ is small enough, then so
is $\|U(0,x)\|_{\widetilde {\textbf{V}}}$. Thus, when
$\|F(0)-w_1\|_{V^*}$ is small enough, by Lemma \ref{l2.6} there
exists a unique solution  $z_0\in \mathscr{O}_{R_2}(0_V)$ for some
$R_2>0$ to
\begin{equation}\label{2.36}
\big(A+A_U(0)+B(0)\big)z_0=F(0)-w_1.
\end{equation}
Since $\|z(0)\|_V\leq c\|z\|_{\mathscr{X}},$ we can choose $R_1$
such that $z(0)\in \mathscr{O}_{R_2}(0_V),$ and we have $z(0)=z_0$.
Therefore, if $\|F\|_{\mathscr{Y}}$ is small enough,
$F(0)-\big(A+A_U(0)+B(0)\big)z_0$ belongs to $H$ and its norm is
small enough, then $z\in\mathscr{X},$ the solution to \eqref{2.35},
is a solution to problem
\begin{equation}\label{2.37}
\left\{\begin{aligned} &
z'(t)+\big(A+A_U(t)+B(t)\big)z(t)=F(t),\\&z(0)= z_0.
\end{aligned}\right.
\end{equation}
 By definitions of
$A, A_U(t), B(t), F,$ the solution $z$ of \eqref{2.37} is also  a
solution to \eqref{2.8} which is equivalent to \eqref{2.3}. Thus,
$e^{k_0t}z+U\in \mathscr{W}$ is a solution to \eqref{2.3}.

 Second, let us prove uniqueness.\\
Let $v_1, v_2$ be two solutions to \eqref{2.3} corresponding to the
same data. Putting $\overline{w}=v_1-v_2$, we have
\begin{equation}\label{2.40}
\begin{aligned}
& \overline{w}\in L_2(0,T;\mathbf{V}),\\& \overline{w}(0)=0,\\&
\langle \overline{w}',u\rangle+
2\nu(\varepsilon(\overline{w}),\varepsilon(u))+\langle(v_1\cdot\nabla)\overline{w},u\rangle+\langle(\overline{w}\cdot\nabla)v_2,u\rangle+2\nu(k(x)\overline{w},u)_{\Gamma_2}\\&
\hspace*{1cm}+2\nu(S\tilde{\overline{w}},\tilde{u})_{\Gamma_3}+2(\alpha(x)\overline{w},u)_{\Gamma_5}+\nu(k(x)\overline{w},u)_{\Gamma_7}=0
\quad
  \mbox{for all}\,\, u\in \textbf{V}.
\end{aligned}
\end{equation}
Putting $w=e^{-k_0t}\overline{w}$, where  $k_0$ is the constant in
Lemma \ref{l2.1}, we get $e^{-k_0t}\overline{w}'=w'+k_0w$. Then, we
have
\begin{equation}\label{2.41}
\begin{aligned}
& w\in L_2(0,T;\mathbf{V}),\\&w(0)=0,\\&
 \langle w',u\rangle+
2\nu(\varepsilon(w),\varepsilon(u))+\langle (v_1\cdot
\nabla)w,u\rangle+\langle (w\cdot
\nabla)v_2,u\rangle+k_0(w,u)\\&\hspace{1cm}+2\nu(k(x)w,u)_{\Gamma_2}+2\nu(S\tilde{w},\tilde{u})_{\Gamma_3}+2(\alpha(x)w,u)_{\Gamma_5}+\nu(k(x)w,u)_{\Gamma_7}
 =0\\
 &\hspace{5cm}\mbox{for all}\,\, u\in \mathbf{V},
\end{aligned}
\end{equation}
which is equivalent to \eqref{2.40}. By Lemma \ref{l2.1},
$$
\begin{aligned}
2\nu(\varepsilon(w),\varepsilon(w))&+k_0(w,w)+2\nu(k(x)w,w)_{\Gamma_2}+2\nu(S\tilde{w},\tilde{w})_{\Gamma_3}\\&+2(\alpha(x)w,w)_{\Gamma_5}+\nu(k(x)w,w)_{\Gamma_7}
\geq \delta \|w\|_\textbf{V}^2.
\end{aligned}
$$
Taking it into account, we can prove $w=0_{\mathscr{X}}$ as in Lemma
3.8 of \cite{k}. Thus, uniqueness of a solution is proved, and we
finished proof of the theorem. $\square$\vspace*{.1cm}

\begin{remark}\label{r2.2}
Let us consider  more precisely the condition that
$F(0)-\big(A+A_U(0)+B(0)\big)z_0$ belongs to $H$ and its norm is
small enough. By \eqref{2.9}$\sim$\eqref{2.12.1} we have
\begin{equation}\label{2.38}
\begin{aligned}
&\big\langle F(0)-\big(A+A_U(0)+B(0)\big)z_0,u\big\rangle=\\
&\Big[-(U'(0,x),u)-2\nu(\varepsilon(U(0,x)),\varepsilon(u))-\langle
(U(0,x)\cdot
\nabla)U(0,x),u\rangle\\&\hspace{.5cm}-2\nu(k(x)U(0,x),u)_{\Gamma_2}-2\nu(S\tilde{U}(0),\tilde{u})_{\Gamma_3}-2(\alpha(x)U(0,x),u)_{\Gamma_5}\\&\hspace{.5cm}-\nu(k(x)U(0,x),u)_{\Gamma_7}+\langle
f(0),u\rangle+\sum_{i=2,4,7}\langle
\phi_i(t),u_n\rangle_{\Gamma_i}+\sum_{i=3,5,6}\langle\phi_i(t),u\rangle_{\Gamma_i}\Big]\\&
-\big[2\nu(\varepsilon(z_0),\varepsilon(u))+2\nu(k(x)z_0,u)_{\Gamma_2}+2\nu(S\tilde{z_0},\tilde{u})_{\Gamma_3}\\&\hspace{4cm}+2(\alpha(x)z_0,u)_{\Gamma_5}+\nu(k(x)z_0,u)_{\Gamma_7}
\big]\\& -\big[ \langle(U(0,x)\cdot
\nabla)z_0,u\rangle+\langle(z_0\cdot\nabla)U(0,x),u\rangle\big]-\langle(z_0,\nabla)z_0,
u\rangle\quad \mbox{for all}\,\, u\in\textbf{ V}.
\end{aligned}
\end{equation}
Taking into account the fact that $U(0,x)+z_0=v_0$, $U'(0,x)\in
\textbf{L}_2(\Omega)$ and its norm is small enough, from
\eqref{2.38} we can see that the condition mentioned above is
equivalent to the condition $\overline{w}_0\in \mathscr{O}_R(0_H)$
for $R>0$ small enough, where $\overline{w}_0$ is defined by
\begin{equation}\label{2.39}
\begin{aligned}
\langle \overline{w}_0, u\rangle&\equiv\langle
f(0),u\rangle+\sum_{i=2,4,7}\langle\phi_i(0,x),u_n\rangle_{\Gamma_i}+\sum_{i=3,5,6}\langle\phi_i(0,x),u\rangle_{\Gamma_i}\\
&-\big[2\nu(\varepsilon(v_0),\varepsilon(u))+2\nu(k(x)v_0,u)_{\Gamma_2}+2\nu(S\tilde{v}_0,\tilde{u})_{\Gamma_3}+2(\alpha(x)v_0,u)_{\Gamma_5}\\&+\nu(k(x)v_0,u)_{\Gamma_7}
+\langle(v_0\cdot \nabla)v_0,u\rangle+k_0(v_0,u)\big]\quad\mbox{for
all}\,\,
 u\in \mathbf{V}.
\end{aligned}
\end{equation}
\end{remark}

\begin{remark}\label{r2.3}
If $\Gamma_i=\varnothing, i=2\sim5, 7,$ then the problem is reduced
to one in \cite{ks} where a local-in-time solution was studied. In
this case $k_0=0$ (cf. Remark \ref{r2.1}), and the condition
\eqref{2.34} is the same with (25) in \cite{ks}. And our condition
for $U$ is also the same with one in \cite{ks}.
\end{remark}\vspace*{.1cm}

\section{Existence of a unique solution to problem II}
\setcounter{equation}{0}

Let
 $\mathbf{V}_1=\{u\in
\mathbf{H}^1(\Omega):\mbox{div}\,u=0,\, u|_{\Gamma_1}=0,\,
u_\tau|_{(\Gamma_2\cup\Gamma_4)}=0,\,
u_n|_{(\Gamma_3\cup\Gamma_5)}=0\}$ and
$\mathbf{V}_{\Gamma2-5}(\Omega)=\{u\in
\mathbf{H}^1(\Omega):\mbox{div}\,u=0,\,
u_\tau|_{(\Gamma_2\cup\Gamma_4)}=0,\,
u_n|_{(\Gamma_3\cup\Gamma_5)}=0\}$. Denote by $H_1$ the completion
of $\mathbf{V}_1$ in the space $\textbf{L}_2(\Omega)$.

By Theorems \ref{t1.1} and \ref{t1.2}, for $v\in
\mathbf{H}^2(\Omega)\cap\mathbf{V}_{\Gamma2-5}(\Omega)$ and $u\in
\mathbf{V}_1$ we have that
\begin{equation}\label{3.1}
\begin{aligned}
-(\Delta v, u)&=(\nabla v, \nabla u)-\left(\f{\p v}{\p
n},u\right)_{\p \Omega}\\& =(\nabla v, \nabla
u)+(k(x)v,u)_{\Gamma_2}-(\mbox{rot}\,v\times n,
u)_{\Gamma_3}+(S\tilde{v},\tilde{u})_{\Gamma_3}\\&
\qquad-(\varepsilon_n(v),u)_{\Gamma_4}
-2(\varepsilon_n(v),u)_{\Gamma_5}-(S\tilde{v},\tilde{u})_{\Gamma_5}-\left(\f{\p
v}{\p n},u\right)_{\Gamma_7}\\& =(\nabla v, \nabla
u)+(k(x)v,u)_{\Gamma_2}-(\mbox{rot}\,v\times n,
u)_{\Gamma_3}+(S\tilde{v},\tilde{u})_{\Gamma_3}\\&
\qquad-(\varepsilon_{nn}(v),u\cdot n)_{\Gamma_4}
-2(\varepsilon_{n\tau}(v),u)_{\Gamma_5}-(S\tilde{v},\tilde{u})|_{\Gamma_5}-\left(\f{\p
v}{\p n},u\right)_{\Gamma_7}.
\end{aligned}
\end{equation}
Also, for $p\in H^1(\Omega)$ and $u\in \mathbf{V}_1$ we get
\begin{equation}\label{3.2}
\begin{aligned}
(\nabla p, u)&=(p, u\cdot n)_{\cup_{i=2}^7\Gamma_i}=(p, u\cdot
n)_{\Gamma_2}+(p, u\cdot n)_{\Gamma_4}+(pn, u)_{\Gamma_7},
\end{aligned}
\end{equation}
where the fact that $u\cdot n|_{\Gamma_3\cup\Gamma_5}=0$ was used.

Let
$$
\begin{aligned}
&\mathscr{X}_1=\{w\in L_2(0,T;\textbf{V}_1); w'\in
L_2(0,T;\textbf{V}_1), w''\in L_2(0,T;\textbf{V}_1^*)\},\\&
\hspace{1cm}\|w\|_{\mathscr{X}_1}=\|w\|_{L_2(0,T;\textbf{V}_1)}+\|w'\|_{L_2(0,T;\textbf{V}_1)}+\|w''\|_{L_2(0,T;\textbf{V}_1^*)},\\
&\mathscr{Y}_1=\{w\in L_2(0,T;\textbf{V}_1^*); w'\in
L_2(0,T;\textbf{V}_1^*)\},\\&
\hspace{1cm}\|w\|_{\mathscr{Y}_1}=\|w\|_{L_2(0,T;\textbf{V}_1^*)}+\|w'\|_{L_2(0,T;\textbf{V}_1^*)}.
\end{aligned}
$$
 Unlike problem I, for problem II we do not require the condition $v_\tau|_{\Gamma_7}=0$,
 and so instead of Assumptions \ref{a2.1} and \ref{a2.2}, we use the
 following assumptions.
\begin{assumption}\label{a3.0}
Assumption \ref{a2.1} holds with $\phi_7, \phi_7' \in L_2(0,T;
\textbf{H}^{-\f{1}{2}}(\Gamma_7))$ instead of $\phi_7, \phi_7' \in
L_2(0,T; H^{-\f{1}{2}}(\Gamma_7))$
\end{assumption}

\begin{assumption}\label{a3.1} There exists
 a function $U\in \mathscr{W}$ such that
$$
\mbox{div}\,U=0,\,U|_{\Gamma_1}=h_1,\, U_\tau|_{\Gamma_2}=0,\,
U_n|_{\Gamma_3}=0,\, U_\tau|_{\Gamma_4}=h_4,\, U_n|_{\Gamma_5}=h_5,
$$
 where
$\mathscr{W}$ is the same as in the previous section. Also,
$U(0,x)-v_0\in \textbf{V}_1$.
\end{assumption}\vspace*{.1cm}

Taking into account \eqref{3.1}, \eqref{3.2}, we get the following
variational formulation for Problem II:
 \begin{formulation}\label{f3.1}
Find $v$ such that
\begin{equation}\label{3.3}
\begin{aligned}
&  v-U\in L_2(0,T;\mathbf{V}_1),\\& v(0)=v_0,\\& \langle
v',u\rangle+ \nu(\nabla v, \nabla
u)+\langle(v\cdot\nabla)v,u\rangle+\nu(k(x)v,u)_{\Gamma_2}\\&\hspace*{2cm}+\nu(S\tilde{v},\tilde{u})_{\Gamma_3}+2(\alpha(x)v,u)_{\Gamma_5}-\nu(S\tilde{v},\tilde{u})_{\Gamma_5}\\&
 \hspace*{1cm}=\langle f,u\rangle+\sum_{i=2,4}\langle \phi_i,u_n\rangle_{\Gamma_i}+\sum_{i=3,5,7}\langle \phi_i,u\rangle_{\Gamma_i}\quad \mbox{for all}\,\, u\in
\mathbf{V}_1.
\end{aligned}
\end{equation}
\end{formulation}\vspace*{.1cm}

Taking into account  Assumption \ref{a3.1} and putting
$v=\overline{z}+U$, we get the following problem which is equivalent
to Formulation \ref{f3.1}:

 Find $\overline{z}$ such that
 \begin{equation}\label{3.4}
\begin{aligned}
& \overline{z}\in L_2(0,T;\mathbf{V}_1),\\&\overline{z}(0)\equiv
v_0-U(0)\in \mathbf{V}_1,\\&
 \langle \overline{z}',u\rangle+
\nu(\nabla \overline{z},\nabla u)+\langle (\overline{z}\cdot
\nabla)\overline{z},u\rangle+\langle (U\cdot
\nabla)\overline{z},u\rangle+\langle (\overline{z}\cdot
\nabla)U,u\rangle\\&\hspace{1cm}+\nu(k(x)\overline{z},u)_{\Gamma_2}+\nu(S\tilde{\overline{z}},\tilde{u})_{\Gamma_3}
+2(\alpha(x)\overline{z},u)_{\Gamma_5}-\nu(S\tilde{\overline{z}},\tilde{u})_{\Gamma_5}\\&
\hspace{.5cm} =-\langle U',u\rangle -\nu(\nabla U,\nabla u)-\langle
(U\cdot
\nabla)U,u\rangle-\nu(k(x)U,u)_{\Gamma_2}-\nu(S\tilde{U},\tilde{u})_{\Gamma_3}\\
&\hspace{1cm}-2(\alpha(x)U,u)_{\Gamma_5}+\nu(S\tilde{U},\tilde{u})_{\Gamma_5}+\langle f,u\rangle
+\sum_{i=2,4}\langle \phi_i,u_n\rangle_{\Gamma_i}+\sum_{i=3,5,7}\langle\phi_i,u\rangle_{\Gamma_i}\\
 &\hspace{6cm}\mbox{for all}\,\, u\in \mathbf{V}_1.
\end{aligned}
\end{equation}

Define an operator $A_{01}:\textbf{V}_1\rightarrow \textbf{V}_1^*$
by
 \begin{equation}\label{3.5}
\begin{aligned}
& \langle A_{01}y,u\rangle=\\& \hspace{.5cm}\nu(\nabla y,\nabla
u)+\nu(k(x)y,u)_{\Gamma_2}+\nu(S\tilde{y},\tilde{u})_{\Gamma_3}+2(\alpha(x)y,u)_{\Gamma_5}-\nu(S\tilde{y},\tilde{u})_{\Gamma_5}\\&
\hspace{5cm}\mbox{for all}\,\, y,u\in \mathbf{V}_1.
\end{aligned}
\end{equation}
By virtue of the same argument used to prove Lemma \ref{l2.1} we get
\begin{lemma}\label{l3.1} $\exists \delta>0$, $\exists k_1\geq 0$:
$$
\langle A_{01}u,u\rangle\geq
\delta\|u\|_{\textbf{V}_1}^2-k_1\|u\|_{H_1}^2\quad\mbox{for all}\,\,
u\in \mathbf{V}_1.
$$
\end{lemma}\vspace{.3cm}

Putting $z=e^{-k_1t}\overline{z}$, where $k_1$ is the constant in
Lemma \ref{l3.1}, and using the fact that
$e^{-k_1t}\overline{z}'=z'+k_1z$, we get the following problem which
is equivalent to \eqref{3.4}:

Find $z$ such that
 \begin{equation}\label{3.6}
\begin{aligned}
& z\in L_2(0,T;\mathbf{V}_1),\\&z(0)=z_0\equiv v_0-U(0)\in
\mathbf{V}_1,\\&
 \langle z'(t),u\rangle+
\nu(\nabla z,\nabla u)+e^{k_1t}\langle (z(t)\cdot
\nabla)z(t),u\rangle+\langle (U(t)\cdot
\nabla)z(t),u\rangle\\&\hspace{1cm}+\langle (z(t)\cdot
\nabla)U(t),u\rangle+k_1(z(t),u)+\nu(k(x)z(t),u)_{\Gamma_2}+\nu(S\tilde{z}(t),\tilde{u})_{\Gamma_3}\\&\hspace{1cm}+2(\alpha(x)z(t),u)_{\Gamma_5}-\nu(S\tilde{z}(t),\tilde{u})_{\Gamma_5}\\&
 =e^{-k_1t}\Big[-(U'(t),u)-\nu(\nabla U,\nabla u)-\langle (U(t)\cdot
\nabla)U(t),u\rangle\\&\hspace{1cm}-\nu(k(x)U(t),u)_{\Gamma_2}
-\nu(S\tilde{U}(t),\tilde{u})_{\Gamma_3}-2(\alpha(x)U(t),u)_{\Gamma_5}
-\nu(S\tilde{U}(t),\tilde{u})_{\Gamma_5}\\&\hspace{1cm}+\langle
f(t),u\rangle +\sum_{i=2,4}\langle \phi_i(t),u_n\rangle_{\Gamma_i}
+\sum_{i=3,5,7}\langle\phi_i(t),u\rangle_{\Gamma_i}\Big]\quad\mbox{for
all}\,\, u\in \mathbf{V}_1.
\end{aligned}
\end{equation}

Define operators $A_1$, $A_{1U}(t)$ by
 \begin{equation}\label{3.7}
\langle A_1v,u\rangle= \langle
A_{01}v,u\rangle+(k_1v,u)\quad\mbox{for all}\,\,
 v, u\in \mathbf{V}_1,
\end{equation}
 \begin{equation}\label{3.8}
\begin{aligned}
&\langle A_{1U}(t)v,u\rangle=  \langle(U(t,x)\cdot
\nabla)v,u\rangle+\langle(v\cdot\nabla)U(t,x),u\rangle\quad\mbox{for
all}\,\,
 v, u\in \mathbf{V}_1,
\end{aligned}
\end{equation}
where $A_{01}$ is one defined in \eqref{3.5}. $U\in \mathscr{W}$
implies $U\in C\left([0,T];\textbf{H}^1(\Omega)\right)$, and such
definitions have meaning. Also, define an operator
$B_1(t):V_1\rightarrow V_1^*$ by
\begin{equation}\label{3.9}
\langle B_1(t)v,u\rangle=e^{k_1t}\langle (v\cdot\nabla)
v,u\rangle\quad\mbox{for all}\,\, v,u\in \mathbf{V}_1.
\end{equation}

 Define an element $F_1\in
\mathscr{Y}_1$ by
 \begin{equation}\label{3.10}
\begin{aligned}
 &\langle
F_1(t),u\rangle=e^{-k_1t}\Big[-\langle U'(t),u\rangle-\nu(\nabla
U(t), \nabla
u)-\langle(U(t)\cdot\nabla)U(t),u\rangle\\&\hspace{0.7cm}-\nu(k(x)U(t),u)_{\Gamma_2}-\nu(S\tilde{U}(t),\tilde{u})_{\Gamma_3}-2(\alpha(x)U(t),u)_{\Gamma_5}+\nu(S\tilde{U}(t),\tilde{u})_{\Gamma_5}\\&
\hspace*{0.7cm}+\langle
f,u\rangle+\sum_{i=2,4}\langle\phi_i,u_n\rangle_{\Gamma_i}+\sum_{i=3,5,7}\langle\phi_i,u\rangle_{\Gamma_i}\Big]\quad
 \mbox{for all}\,\, u\in \mathbf{V}_1.
\end{aligned}
\end{equation}

Now,  in the same way as Theorem \ref{t2.1} we can prove the
following theorem which is one of the main results of this paper.
\begin{theorem}\label{t3.1}
Suppose that Assumptions \ref{a3.0} and \ref{a3.1} hold.
 Assume that $\|U\|_{\mathscr{W}}$ and the norms of
 $f, f', \phi_i, \phi_i'$ in the spaces they belong to are small enough.

If
\begin{equation}\label{3.11}
\begin{aligned}
w_1\equiv F_1(0)-(A_1+A_{1U}(0)+B_1(0))z_0\in H_1,
\end{aligned}
\end{equation}
where $z_0=v_0-U(0,\cdot),$ and $\|w_1\|_{H_1}$ is small enough,
then in the space $\mathscr{W}$ there exists a unique solution to
\eqref{3.3}.
\end{theorem}

\begin{remark}\label{r3.1}
By the same argument as Remark \ref{r2.2}, we can see that the
condition \eqref{3.11} is equivalent to the condition
$\overline{w}_0\in H_1,$ where $\overline{w}_0\in \mathbf{V}_1^*$ is
defined by
 \begin{equation}\label{3.11.1}
\begin{aligned}
\langle \overline{w}_1, u\rangle&=\langle
f(0),u\rangle+\sum_{i=2,4}\langle\phi_i(0,x),u_n\rangle_{\Gamma_i}+\sum_{i=3,5,7}\langle\phi_i(0,x),u\rangle_{\Gamma_i}\\
&-\big[\nu(\nabla v_0,\nabla
u)+\nu(k(x)v_0,u)_{\Gamma_2}+\nu(S\tilde{v}_0,\tilde{u})_{\Gamma_3}+2(\alpha(x)v_0,u)_{\Gamma_5}\\&-\nu(S\tilde{v}_0,\tilde{u})_{\Gamma_5}
+\langle(v_0\cdot \nabla)v_0,u\rangle+k_1(v_0,u)\big]\quad\mbox{for
all}\,\,
 u\in \mathbf{V}_1,
\end{aligned}
\end{equation}
with $k_1$ in Lemma \ref{l3.1}.
\end{remark}

\begin{remark}\label{r3.2}
If $U\equiv 0$ and $\Gamma_i=\varnothing, i=2\sim 5,$ then problem
\eqref{3.3} is reduced to one in \cite{k}.  In this case $k_1=0.$
(cf. Remark \ref{r2.1}). If $v_0\in \textbf{H}^{l/2}(\Omega)$, then
$(v_0\cdot \nabla)v_0\in \textbf{L}_2(\Omega)$. Thus, the condition
above for $\overline{w}_1$ being in $H_1$ is the same with one of
conditions of Theorems 3.5$\sim$3.8 in \cite{k}, but we do not
demand $v_0\in \textbf{H}^{r_0}(\Omega),\, r_0>\f{l}{2}$.
\end{remark}

\section{Existence of a unique solution for perturbed data}
\setcounter{equation}{0}

In \cite{k} it is proved that if a solution satisfying smoothness
and a compatibility condition is given, then there exists a unique
solution for small perturbed data satisfying the compatibility
condition. In this section we get such results for the Problems I
and II. In our results the conditions for a given solution is
essentially the same with one in \cite{k}, but the smoothness
condition for the initial functions in the compatibility condition
for small perturbed data is weaker than one in \cite{k}(cf. Remark
\ref{r4.1}).

Let $\widetilde{\textbf{V}}^{r_0}=\{u\in
\mathbf{H}^{r_0}(\Omega):\mbox{div}\,u=0\}$, $r_0>l/2$, and
$$
\begin{aligned}
&\overline{\mathscr{W}}=\left\{w\in L_2(0,T;\widetilde{\textbf{V}});
w'\in L_2(0,T;\widetilde{\textbf{V}}), w''\in L_2(0,T;{\widetilde
{\textbf{V}}}^*),w(0)\in \widetilde{\textbf{V}}^{r_0}\right\},\\&
\hspace{1cm}\|w\|_{\mathscr{W}}=\|w\|_{L_2(0,T;\widetilde{\textbf{V}})}+\|w'\|_{L_2(0,T;\widetilde{\textbf{V}})}+
\|w''\|_{L_2(0,T;{\widetilde
{\textbf{V}}}^*)}+\|w(0)\|_{\widetilde{\textbf{V}}^{r_0}}.
\end{aligned}
$$

Let us consider Problem I.

Let $W(x,t)\in \overline{\mathscr{W}}$ be a given solution to
Problem I. Let $v$ be the solution for the data perturbed except
$h_i$ and put $v=\overline{z}+W$. Then, we get a problem for
$\overline{z}$:

 Find $\overline{z}$ such that
 \begin{equation}\label{4.1}
\begin{aligned}
& \overline{z}\in L_2(0,T;\mathbf{V}),\\&\overline{z}(0)=z_0\equiv
v_0-W(0,x)\in \mathbf{V},\\&
 \langle \overline{z}',u\rangle+
2\nu(\varepsilon(\overline{z}),\varepsilon(u))+\langle
(\overline{z}\cdot \nabla)\overline{z},u\rangle+\langle (W\cdot
\nabla)\overline{z},u\rangle+\langle (\overline{z}\cdot
\nabla)W,u\rangle\\&
\hspace{1cm}+2\nu(k(x)\overline{z},u)_{\Gamma_2}+2\nu(S\tilde{\overline{z}},\tilde{u})_{\Gamma_3}+2(\alpha(x)\overline{z},u)_{\Gamma_5}+\nu(k(x)\overline{z},u)_{\Gamma_7}\\&
 =\langle f,u\rangle+\sum_{i=2,4,7}\langle \phi_i,u_n\rangle_{\Gamma_i}+\sum_{i=3,5,6}\langle\phi_i,u\rangle_{\Gamma_i}\quad\mbox{for all}\,\, u\in \mathbf{V},
\end{aligned}
\end{equation}
where $z_0, f, \phi_i$ are perturbations of corresponding data.

\begin{remark}\label{r4.0}
Proofs of this section are similar to one in Section 3. Main
difference is that we do not assume smallness of $W(0,x)$ unlike
$U(0,x)$ in Section 3.
\end{remark}\vspace{.1cm}

Define an operator $A_{02}: \textbf{V}\rightarrow \textbf{V}^*$ by
\begin{equation}\label{4.2}
\begin{aligned}
 \langle A_{02}y,u\rangle=&
2\nu(\varepsilon(y),\varepsilon(u))+2\nu(k(x)y,u)_{\Gamma_2}+2\nu(S\tilde{y},\tilde{u})_{\Gamma_3}\\&+2(\alpha(x)y,u)_{\Gamma_5}+\nu(k(x)y,u)_{\Gamma_7}+\langle
(W(0,x)\cdot \nabla)y,u\rangle\\&+\langle (y\cdot
\nabla)W(0,x),u\rangle\qquad \mbox{for all}\,\, y,u\in \mathbf{V}.
\end{aligned}
\end{equation}
\begin{lemma}\label{l4.1} There exists $\delta>0$ and $k_2\geq 0$ such that
$$
\langle A_{02}u,u\rangle\geq
\delta\|u\|_\textbf{\textbf{V}}^2-k_2\|u\|_H^2\quad\mbox{for
all}\,\, u\in \mathbf{V}.
$$
\end{lemma}
\emph{Proof}\,\, By Korn's inequality
\begin{equation}\label{4.3}
2\nu(\varepsilon(u),\varepsilon(u))\geq \beta
\|u\|^2_\mathbf{V}\quad\exists\beta>0, \mbox{for all}\,\,
 u\in \mathbf{V}.
\end{equation}
By Remark \ref{r1.1}, there exists  a constant $M$ such that
\[
\begin{aligned}
\|S(x)\|_\infty,\, \|k(x)\|_\infty, \|\alpha(x)\|_\infty \leq M.
\end{aligned}
\]
Then, there exists a constant $c_0$ (depend on $\beta$) such that
\begin{equation}\label{4.4}
\begin{aligned}
&\left|2\nu(k(x)u,u)_{\Gamma_2}
+2\nu(S\tilde{u},\tilde{u})_{\Gamma_3}+\nu(k(x)u,u)_{\Gamma_7}+2(\alpha(x)u,u)_{\Gamma_5}\right|\\&
\hspace{2cm}\leq \f{\beta}{4}\|
u\|^2_{\mathbf{H}^1(\Omega)}+c_0\|u\|_{H}^2\quad \mbox{for all}\,\,
u\in \mathbf{V}
\end{aligned}
\end{equation}
(cf. Theorem 1.6.6 in \cite{bs}).

Let us estimate
$\langle(W(0,x)\cdot \nabla)u,u\rangle+\langle (u\cdot
\nabla)W(0,x),u\rangle$. Since $W(0,x)\in C(\overline{\Omega})$,
\begin{equation}\label{4.5}
\big|\big\langle(W(0,x)\cdot \nabla)u,u\big\rangle\big|\leq
\f{\beta}{8}\| u\|^2_{\mathbf{H}^1(\Omega)}+c_1\|u\|_{H}^2.
\end{equation}
Taking $\text{div}\, u=0$ into account, we get
$$
\begin{aligned}
\big\langle (u\cdot \nabla)W(0,x),u\big\rangle&=\sum_j\int_\Omega
u_j\f{\p W(0,x)}{\p x_j}u\, dx\\&=\int_{\p\Omega}(u\cdot
n)(W(0,x)\cdot u)\,d\Gamma-\sum_j\int_\Omega u_j\f{\p u}{\p
x_j}W(0,x)\,dx.
\end{aligned}
$$
On the right hand side of the formula above estimating the first
term as in \eqref{4.4} and applying H\"{o}lder inequality in the
second term, we have
\begin{equation}\label{4.6}
\big|\big\langle (u\cdot \nabla)W(0,x),u\big\rangle\big|\leq
\f{\beta}{8}\| u\|^2_{\mathbf{H}^1(\Omega)}+c_2\|u\|_{H}^2.
\end{equation}
Putting $\delta=\f{\beta}{2}$, $k_2=c_0+c_1+c_2$, from
\eqref{4.3}-\eqref{4.6} we get the asserted conclusion.
$\square$\vspace*{.3cm}

Put $z=e^{-k_2t}\overline{z}$, where $k_2$ is a constant in Lemma
\ref{l4.1}. Then, $e^{-k_2t}\overline{z}'=z'+k_2z$ and we have the
following problem which is equivalent to \eqref{4.1}.

Find $z$ such that
\begin{equation}\label{4.7}
\begin{aligned}
& z\in L_2(0,T;\mathbf{V}),\\&z(0)=z_0=v_0-W(0)\in \mathbf{V},\\&
 \langle z'(t),u\rangle+
2\nu(\varepsilon(z)(t),\varepsilon(u))+e^{k_2t}\langle (z(t)\cdot
\nabla)z(t),u\rangle+\langle (W(t)\cdot
\nabla)z(t),u\rangle\\&\hspace{1cm}+\langle (z(t)\cdot
\nabla)W(t),u\rangle+k_2(z(t),u)+2\nu(k(x)z(t),u)_{\Gamma_2}+2\nu(S\tilde{z}(t),\tilde{u})_{\Gamma_3}\\&\hspace{1cm}+2(\alpha(x)z(t),u)_{\Gamma_5}+\nu(k(x)z(t),u)_{\Gamma_7}\\&
 =e^{-k_2t}\Big[\langle f(t),u\rangle+\sum_{i=2,4,7}\langle \phi_i(t),u_n\rangle_{\Gamma_i}+\sum_{i=3,5,6}\langle\phi_i(t),u\rangle_{\Gamma_i}\Big]\quad\mbox{for all}\,\, u\in \mathbf{V}.
\end{aligned}
\end{equation}

Define operators $A_2, A_W(t):\textbf{V}\rightarrow \textbf{V}^*$ by
 \begin{equation}\label{4.8}
\begin{aligned}
 \langle A_2y,u\rangle=&
2\nu(\varepsilon(y),\varepsilon(u))+2\nu(k(x)y,u)_{\Gamma_2}+2\nu(S\tilde{y},\tilde{u})_{\Gamma_3}\\&+2(\alpha(x)y,u)_{\Gamma_5}+\nu(k(x)y,u)_{\Gamma_7}+k_2(y,u)\quad
\mbox{for all}\,\, y,u\in \mathbf{V},
\end{aligned}
\end{equation}
 \begin{equation}\label{4.9}
\begin{aligned}
&\langle A_W(t)v,u\rangle=\langle(W(t,x)\cdot
\nabla)v,u\rangle+\langle(v\cdot\nabla)W(t,x),u\rangle\quad\mbox{for
all}\,\,
 v, u\in \mathbf{V},
\end{aligned}
\end{equation}
where $k_2$ is a constant in Lemma \ref{l4.1}. $W\in
\overline{\mathscr{W}}$ implies $W\in
C\left([0,T];\textbf{H}^1(\Omega)\right)$, and such definitions are
well.

In proof of Lemma \ref{l4.1} it is clear that
 \begin{equation}\label{4.10}
\left\langle A_2u,u\right\rangle \geq \f{3\beta}{4}
\|u\|^2_\textbf{V}.
\end{equation}
Also, by Lemma \ref{l4.1}
\begin{equation}\label{4.11} \left\langle
\big(A_2+A_W(0)\big)u,u\right\rangle \geq \f{\beta}{4
}\|u\|^2_\textbf{V}.
\end{equation}
Define an operator $B_2(t):V\rightarrow V^*$ by
\begin{equation}\label{4.12}
\langle B_2(t)v,u\rangle=e^{k_2t}\langle (v\cdot\nabla)
v,u\rangle\quad\mbox{for all}\,\, v,u\in \mathbf{V}.
\end{equation}
Define operators $L_2, {\widetilde A}_W, L_{2W}, {\widetilde
B}_2:\mathscr{X}\rightarrow\mathscr{Y}$,
$C_2:\mathscr{X}\times\mathscr{X}\rightarrow\mathscr{Y}$ and an
element $F_2\in \mathscr{Y}$ by
 \begin{equation}\label{4.13}
\begin{aligned}
&\langle (L_2z)(t),u\rangle=\langle z'(t),u\rangle+\langle
A_2z(t),u\rangle \quad\mbox{for all}\,\,
z\in \mathscr{X}, \mbox{for all}\,\, u\in \mathbf{V},\\
&\langle ({\widetilde A}_Wz)(t),u\rangle=\langle A_W(t)z(t),u\rangle
\quad\mbox{for all}\,\, z\in \mathscr{X}, \mbox{for all}\,\, u\in
\mathbf{V},\\& \langle(L_{2W}z)(t),u\rangle=\langle
z'(t),u\rangle+\big\langle \big(A_2+A_W(t)\big)z(t),u\big\rangle
\quad\mbox{for all}\,\,
z\in \mathscr{X}, \mbox{for all}\,\, u\in \mathbf{V},\\
&\langle ({\widetilde B}_2z)(t),u\rangle=\langle B_2(t)z(t), u
\rangle\quad\mbox{for all}\,\,
z\in \mathscr{X}, \mbox{for all}\,\, u\in \mathbf{V},\\
&\langle C_2(w,z)(t),u\rangle=e^{k_2t}\langle (w\cdot\nabla)
z,u\rangle+e^{k_2t}\langle (z\cdot\nabla)w,u\rangle\quad\mbox{for
all}\,\,
w,z\in \mathscr{X}, \mbox{for all}\,\, u\in \mathbf{V},\\
 &\langle
(F_2)(t),u\rangle=e^{-k_2t}\Big[\langle
f(t),u\rangle+\sum_{i=2,4,7}\langle
\phi_i(t),u_n\rangle_{\Gamma_i}+\sum_{i=3,5,6}\langle\phi_i(t),u\rangle_{\Gamma_i}\Big]\quad\mbox{for
all}\,\, u\in \mathbf{V}.
\end{aligned}
\end{equation}

By the argument as Lemma \ref{l2.2} we get
 \begin{lemma}\label{l4.2}  $C_2$ is a bilinear continuous operator such that
 $\mathscr{X}\times \mathscr{X}\rightarrow\mathscr{Y}$.
 Under Assumption \ref{a2.1} $ {\widetilde A}_W$ is a linear continuous operator such that
  $\mathscr{X}\rightarrow\mathscr{Y}$ and $F_2\in \mathscr{Y}$.
\end{lemma}

Using \eqref{4.10} instead of \eqref{2.18}, as Lemma \ref{l2.3} we
get
\begin{lemma}\label{l4.3} The operator $\overline{L}_2$ defined by
$\overline{L}_2z=(z'(0), L_2z)$ for $z\in \mathscr{X}$ is a linear
continuous one-to-one operator from $\mathscr{X}$ onto $H\times
\mathscr{Y}$.
\end{lemma}

Now, using \eqref{4.11} without assuming the fact that
$\|W(0,x)\|_{\widetilde {\textbf{V}}}$ is small enough, as Lemma
\ref{l2.4} we prove the following
\begin{lemma}\label{l4.4}
 The operator $\overline{L}_{2W}$
defined by $\overline{L}_{2W}z=(z'(0), L_{2W}z)$ for $z\in
\mathscr{X}$ is a linear continuous one-to-one operator from
$\mathscr{X}$ onto $H\times \mathscr{Y}$.
\end{lemma}

\emph{Proof}\,\,  As Lemma 3.5 in \cite{k} it is proved that the
operator ${\widetilde A}_W\in (\mathscr{X}\rightarrow\mathscr{Y})$
is compact. Thus, $\overline{A}_W\in (\mathscr{X}\rightarrow
H\times\mathscr{Y})$ defined by $\overline{A}_Wz=\{0_H,
\widetilde{A}_Wz\}$ is also compact. Since
$\overline{L}_{2W}=\overline{L}_2+\overline{A}_W$, in order to get
the asserted conclusion by virtue of Theorem 3.4 in \cite{k}  it is
enough to prove that $\overline{L}_{2W}$ is one-to-one from
$\mathscr{X}$ into $H\times \mathscr{Y}$.

To prove the last fact it is enough to show that the inverse image
of $(0_H, 0_{\mathscr{Y}})$ by $\overline{L}_{2W}$ is
$0_{\mathscr{X}}$. It is easy to verify that
\begin{equation}\label{4.14}
 \big|\big\langle \big(A_2+A_W(0)\big)v,u\big\rangle\big|\leq c\|v\|_\mathbf{V}\cdot \|u\|_\mathbf{V}
\quad\mbox{for all}\,\,
 v,u\in \mathbf{V}.
\end{equation}
By \eqref{4.11}, \eqref{4.14} for any $q\in \textbf{V}^*$ there
exists a unique solution  $y\in V$ to
\begin{equation}\label{4.15}
(A_2+A_W(0))y=q.
\end{equation}

Let $z\in\mathscr{X}$  be the inverse image of $(0_H,0_\mathscr{Y})$
by $\overline{L}$. Then, $z'(0)=0_H$, and putting $t=0$ from the
third one in \eqref{4.13} we get
$$
\big\langle \big(A_2+A_W(0)\big)z(0),u\big\rangle=0 \quad\mbox{for
all}\,\,
 u\in \mathbf{V},
$$
where $z(0)=z(0,x)$. This means that $z(0)$ is the unique solution
to \eqref{4.15} with $q=0_{\textbf{V}^*}$, i.e.,
$z(0)=0_\textbf{V}$. Therefore, $z\in \mathscr{X}$ satisfies
\begin{equation}\label{4.16}
\left\{\begin{aligned}
&z'(t)+\big(A_2+A_W(t)\big)z(t)=0,\\
 &z(0)=0_\textbf{V}.
 \end{aligned}\right.
\end{equation}
Now, using \eqref{4.16} and Gronwall's inequality, as in Lemma 3.8
of \cite{k} we can prove $z=0_{\mathscr{X}}$.
 It is finished to prove the Lemma.
 $\square$\vspace*{.1cm}\vspace*{.1cm}

By the argument as Lemma \ref{l2.5} we get
\begin{lemma}\label{l4.5} The operator $T_2$ defined by $T_2z=\left(z'(0), (L_{2W}+{\widetilde B}_2)z\right)$ for $z\in \mathscr{X}$ is continuously
differentiable, $T_20_\mathscr{X}=\left(0_H,0_\mathscr{Y}\right)$
and the Frechet derivative of $T_2$ at $0_\mathscr{X}$ is
$\overline{L}_{2W}$.
\end{lemma}

Let us consider the following problem
\begin{equation}\label{4.17}
\big(A_2+A_W(0)+B_2(0)\big)u=q.
\end{equation}

Now, using \eqref{4.11} without assuming the fact that
$\|W(0,x)\|_{\widetilde {\textbf{V}}}$ is small enough, as Lemma
\ref{l2.6} we can prove
\begin{lemma}\label{l4.6} If the norm of $q\in V^*$ is small enough, then there exists a
unique solution to \eqref{4.17} in some
$\mathscr{O}_M(0_\mathbf{V})$.
\end{lemma}

Using Lemmas \ref{l4.2}$\sim$\ref{l4.5}, Proposition \ref{p2.1}, in
the same way as Theorem \ref{t2.1} we get
\begin{theorem}\label{t4.1}
Suppose that Assumptions \ref{a2.1} holds and the norms of
 $f, f', \phi_i,$ $ \phi_i'$ in the spaces they belong to are small enough.

If
\begin{equation}\label{4.18}
\begin{aligned}
w_2\equiv F_2(0)-(A_2+A_{2W}(0)+B_2(0))z_0\in H,
\end{aligned}
\end{equation}
where $z_0=v_0-U(0,\cdot),$ and $\|w_2\|_{H}$ is small enough, then
there exists a unique solution to \eqref{4.1} in the space
$\mathscr{W}.$
\end{theorem}

\begin{remark}\label{r4.1}
By the same argument as Remark \ref{r2.2}, we can see that the
condition \eqref{4.18} is equivalent to the condition
$\overline{w}_2\in H_1,$ where $\overline{w}_2\in \mathbf{V}_1^*$ is
defined  by
\begin{equation}\label{4.18-0}
\begin{aligned}
\langle \overline{w}_2, u\rangle&=\langle
f(0),u\rangle+\sum_{i=2,4}\langle\phi_i(0,x),u_n\rangle_{\Gamma_i}+\sum_{i=3,5,7}\langle\phi_i(0,x),u\rangle_{\Gamma_i}\\
&-\big[2\nu(\varepsilon(z_0),\varepsilon(u))+2\nu(k(x)z_0,u)_{\Gamma_2}+2\nu(S\tilde{z}_0,\tilde{u})_{\Gamma_3}\\&
+2(\alpha(x)z_0,u)_{\Gamma_5}+\nu(k(x)z_0,u)_{\Gamma_7}+\langle
(W(0,x)\cdot \nabla)z_0,u\rangle\\&+\langle (z_0\cdot
\nabla)W(0,x),u\rangle+k_2(z_0,u)+\langle(z_0\cdot
\nabla)z_0,u\rangle\big]\quad\mbox{for all}\,\,
 u\in \mathbf{V}
\end{aligned}
\end{equation}
with $k_2$ in Lemma \ref{l4.1}.
\end{remark}

Let us consider Problem II.

Let $W(x,t)\in \overline{\mathscr{W}}$ be given solution to Problem
II. Let $v$ be the solution for the data perturbed except $h_i$ and
put $v=\overline{z}+W$. Then, we get a problem for $\overline{z}$:

Find $\overline{z}$ such that
\begin{equation}\label{4.19}
\begin{aligned}
& \overline{z}\in L_2(0,T;\mathbf{V}_1),\\&\overline{z}(0)=z_0\equiv
v_0-W(0,x)\in \mathbf{V}_1,\\& \langle \overline{z}',u\rangle+
\nu(\nabla \overline{z},\nabla u)+\langle (\overline{z}\cdot
\nabla)\overline{z},u\rangle+\langle (W\cdot
\nabla)\overline{z},u\rangle+\langle (\overline{z}\cdot
\nabla)W,u\rangle\\&\hspace{1cm}+\nu(k(x)\overline{z},u)_{\Gamma_2}+\nu(S\tilde{\overline{z}},\tilde{u})_{\Gamma_3}+2(\alpha(x)\overline{z},u)_{\Gamma_5}-\nu(S\tilde{\overline{z}},\tilde{u})_{\Gamma_5}\\&
\hspace{.5cm} =\langle f,u\rangle+\sum_{i=2,4}\langle
\phi_i,u_n\rangle_{\Gamma_i}+\sum_{i=3,5,7}\langle\phi_i,u\rangle_{\Gamma_i}
\quad\mbox{for all}\,\, u\in \mathbf{V}_1,
\end{aligned}
\end{equation}
where $z_0, f, \phi_i$ are perturbations of corresponding data.

By the same argument as Theorem \ref{t4.1} we have
\begin{theorem}\label{t4.2}
Suppose that Assumptions \ref{a2.1} holds and the norms of
 $f, f',$ $ \phi_i, \phi_i'$ in the spaces they belong to are small enough.
Define an element $w_3\in \mathbf{V}_1^*$ by
\begin{equation}\label{4.20}
\begin{aligned}
\langle w_3, u\rangle&=\langle
f(0),u\rangle+\sum_{i=2,4}\langle\phi_i(0,x),u_n\rangle_{\Gamma_i}+\sum_{i=3,5,7}\langle\phi_i(0,x),u\rangle_{\Gamma_i}\\
&-\big[\nu(\nabla z_0,\nabla
u)+\nu(k(x)z_0,u)_{\Gamma_2}+\nu(S\tilde{z}_0,\tilde{u})_{\Gamma_3}+2(\alpha(x)z_0,u)_{\Gamma_5}
\\&-\nu(S\tilde{z}_0,\tilde{u})_{\Gamma_5}+\langle (W(0,x)\cdot \nabla)z_0,u\rangle+\langle
(z_0\cdot \nabla)W(0,x),u\rangle\\&+k_3(z_0,u)+\langle(z_0\cdot
\nabla)z_0,u\rangle\big]\hspace{.4cm}\mbox{for all}\,\,
 u\in \mathbf{V}_1,
\end{aligned}
\end{equation}
where $k_3$ is a constant determined as in Lemma \ref{l4.1}.

If $w_3 \in \mathscr{O}_R(0_{H_1})$ for $R>0$ small enough, then
there exists a unique solution to \eqref{4.19} in the space
$\mathscr{W}$.
\end{theorem}

\begin{remark}\label{r4.2}
If $\Gamma_i=\varnothing, i=2\sim5,$ then problem \eqref{4.19} is
reduced to one in \cite{k}. If $z_0\in \textbf{H}^{l/2}(\Omega)$,
then $(z_0\cdot \nabla)z_0, (W(0,x)\cdot \nabla)z_0, (z_0\cdot
\nabla)W(0,x)\in \textbf{L}_2(\Omega)$ and $k_3z_0\in
\textbf{L}_2(\Omega)$. Thus, the last 4 terms in the right hand side
of \eqref{4.20} do not give any effect to the condition for $w_3$
being in $H_1$, and so the conditions in the Theorem \ref{t4.2} is
the same with one of conditions of Theorems 3.5$\sim$3.8 in
\cite{k}. Thus, Theorem \ref{t4.2} guarantees existence of a unique
solution under a condition weaker that one in \cite{k}.

Note that putting $W(t,x)\equiv 0$ in Theorems \ref{t4.1} and
\ref{t4.2}, we can not get Theorems \ref{t2.1} and \ref{t3.1}, since
there $h_i\neq 0$.
\end{remark}\vspace*{.1cm}

{\bf Acknowledgment:}\,\,The authors are grateful to the anonymous
referee for his or her valuable comments.

\end{document}